\documentclass[11pt,letterpaper]{article}
\usepackage[margin=1in]{geometry}
\usepackage[T1]{fontenc}
\usepackage[utf8]{inputenc}
\usepackage{amsmath,amssymb,amsthm,mathtools}
\usepackage{booktabs,array,enumitem,microtype}
\usepackage[table]{xcolor}
\usepackage{xspace,threeparttable,needspace}
\usepackage{algorithm,algpseudocode,etoc}
\usepackage[most]{tcolorbox}
\usepackage[round,authoryear]{natbib}
\setcitestyle{citesep={,}}
\definecolor{citationred}{HTML}{9B3348}
\usepackage[colorlinks=true,citecolor=citationred,linkcolor=black,urlcolor=black]{hyperref}
\newtheorem{theorem}{Theorem}[section]
\newtheorem{lemma}[theorem]{Lemma}
\newtheorem{proposition}[theorem]{Proposition}
\newtheorem{corollary}[theorem]{Corollary}
\definecolor{statementbackground}{HTML}{F3F1EB}
\tcbset{statementshade/.style={enhanced,breakable,colback=statementbackground,
 colframe=blue!15!black,boxrule=0pt,arc=1pt,left=6pt,right=6pt,
 top=4pt,bottom=4pt,before skip=8pt,after skip=8pt}}
\tcolorboxenvironment{theorem}{statementshade}
\tcolorboxenvironment{lemma}{statementshade}
\tcolorboxenvironment{proposition}{statementshade}
\tcolorboxenvironment{corollary}{statementshade}
\theoremstyle{definition}
\newtheorem{assumption}[theorem]{Assumption}
\newtheorem{example}[theorem]{Example}
\tcolorboxenvironment{example}{statementshade,coltext=black}
\theoremstyle{remark}
\newtheorem{remark}[theorem]{Remark}
\newcommand{\X}{\mathcal{X}}
\newcommand{\Xs}{\mathcal{X}_\star}
\newcommand{\R}{\mathbb{R}}
\newcommand{\B}{\mathbb{B}}
\newcommand{\dist}{\operatorname{dist}}
\newcommand{\Proj}{\operatorname{Proj}}
\newcommand{\argmin}{\operatorname*{arg\,min}}
\newcommand{\algname}[1]{{\small\sffamily #1}\xspace}
\newcommand{\algnametiny}[1]{{\sffamily #1}\xspace}

\newcommand{\pos}[1]{\left[#1\right]_+}
\newcommand{\ceilp}[1]{\left\lceil\pos{#1}\right\rceil}
\newcommand{\inner}[2]{\langle #1,#2\rangle}
\newenvironment{relchain}{\begin{array}{@{}>{\displaystyle}r@{\;}>{\displaystyle}c@{\;}>{\displaystyle}l@{}}}{\end{array}}
\newenvironment{revision}{\begingroup\hypersetup{citecolor=blue,linkcolor=blue,urlcolor=blue}}{\par\endgroup}
\DeclareRobustCommand{\rev}[1]{{\hypersetup{citecolor=blue,linkcolor=blue,urlcolor=blue}#1}}

\DeclareRobustCommand{\auditrev}[1]{#1}
\DeclareRobustCommand{\simplerev}[1]{#1}
\DeclareRobustCommand{\termrev}[1]{\textcolor{black}{#1}}
\setlist{leftmargin=*,topsep=4pt,itemsep=4pt}
\title{\bfseries An Adaptive Short-Step Variant of BroxGD\\with Linear Convergence under Quadratic Growth}
\author{Jan Niklas Kolzenburg\textsuperscript{1,2}\thanks{Jan Niklas Kolzenburg was affiliated with Frankfurt School of Finance \& Management. This paper was written during his internship in the KAUST Optimization \& Machine Learning Lab. All results were obtained by early September 2026. A September 2026 revision of the \algname{BroxGD} paper \citep{broxgd2026v2} subsequently proposed several adaptive variants of \algname{BroxGD}; our results preceded that development.}\qquad
Peter Richt{\'a}rik\textsuperscript{1}\\[5pt]
\normalsize\textsuperscript{1}King Abdullah University of Science and Technology (KAUST)\\
\normalsize\textsuperscript{2}Frankfurt School of Finance \& Management}
\date{7 November 2026}

\begin{document}
\etocdepthtag.toc{main}
{\let\newpage\relax\maketitle}
\begin{abstract}
Broximal Gradient Descent  (\algname{BroxGD}), recently analyzed by \citet{locallmo}, is a projection-free sister method to Projected Gradient Descent (\algname{ProjGD}) for solving problems of the form
$\min_{x \in \X} f(x),$
where $f:\R^d \to \R$, and $\X \subseteq \R^d$ is a nonempty, closed and convex set. The method performs {\em local} linear minimization  of the form
\[x_{k+1} \in \argmin_{z\in\X\cap\B(x_k,t_k)}\inner{\nabla f(x_k)}{z},\]
where $\B(x,t)$ is the Euclidean ball centered at $x$ of radius $t$, and $\{t_k\}$ are positive radii chosen to facilitate convergence. Just like \algname{ProjGD}, \algname{BroxGD} enjoys a linear convergence guarantee in the smooth strongly convex setting, {with the same condition-number dependence on the ratio of the smoothness and strong convexity constants}. {In contrast, vanilla Frank--Wolfe, which relies on global linear minimization over $\X$, does not generally converge linearly in this regime, even on polytopes; variants such as away-step Frank--Wolfe obtain linear convergence under additional geometric assumptions.}  However, the radius rule for which \citet{locallmo} prove  their result  requires unavailable information---the distance of the current iterate to the optimal solution---and is therefore impractical.   {We address this limitation by developing an adaptive short-step variant of \algname{BroxGD}} which converges linearly for smooth {convex} functions under quadratic growth (which generalizes strong convexity), and which does not require the knowledge of any unavailable information.  Our short-step method has the simple form 
  \[z_{k} {\in} \argmin_{z\in\X\cap\B(x_k,t_k)}\inner{\nabla f(x_k)}{z}, \quad x_{k+1} = \Proj_{S_k} \left(x_k-\nabla f(x_k)/L\right),\]
  where $L$ is the smoothness constant of $f$, $S_k$ is the line segment connecting $x_k$ and $z_{k}$,  and the radius is updated via  the simple adaptive rule
$t_{k+1}=\max\left\{
\|x_{k+1}-x_k\|,\,
t_k/2
\right\}.$
{Projection onto $S_k$ requires only scalar clipping, not projection onto $\X$.}

%
\begingroup\color{black}
Assuming exact access to this local linear minimization oracle, each update uses one oracle call. For any positive initial radius, our explicit bound counts all oracle calls needed to reach objective error at most $\varepsilon$ and yields $O\!\left((L/\alpha)\log(1/\varepsilon)\right)$ for fixed initialization, where $\alpha=\min\left\{\alpha_0,L\right\}$ and $\alpha_0>0$ is the quadratic-growth constant. The algorithm does not require $\alpha_0$.
\par\endgroup
\end{abstract}

\clearpage
\begingroup
\etocsettagdepth{main}{subsection}
\etocsettagdepth{appendix}{subsection}
\renewcommand{\rev}[1]{#1}
\tableofcontents
\endgroup
\clearpage

\section{Introduction}
We consider the optimization problem
\begin{equation}\label{eq:problem}
 \min_{x\in\X}f(x),
\end{equation}
where $f$ is smooth and convex and $\X\subseteq\R^d$ is nonempty, closed, and convex. Our linear-convergence results assume quadratic growth and therefore also apply to strongly convex objectives. Projected gradient descent (\algname{ProjGD}) is a standard approach to this problem. It projects each gradient step onto $\X$ \citep{beck2017firstorder}. When this projection is expensive, methods that maintain feasibility through a different oracle offer an alternative. Frank--Wolfe (\algname{FW}) minimizes the linearized objective over $\X$ and moves along the segment from the current iterate toward the resulting minimizer \citep{FW,jaggi13}. Linear-convergence guarantees for variants using this global oracle typically rely on additional structure in the objective or feasible set, as in away-step methods on polytopes \citep{lacostejulien2015}.

\begingroup\color{black}
Broximal Gradient Descent (\algname{BroxGD}) replaces projection by local linear minimization \citep{locallmo}. For a linear objective vector $g\in\R^d$, a center $x\in\X$, and a radius $t>0$, its oracle solves
\begin{equation}\label{eq:oracle}
 \argmin_{z\in\X\cap\B(x,t)}\inner{g}{z},
 \qquad \B(x,t)=\left\{z:\|z-x\|\le t\right\}.
\end{equation}
Algorithm~\ref{alg:main} uses $g=\nabla f(x_k)$. The ball makes this subproblem compact even when $\X$ is unbounded. The same Euclidean subproblem appears in the earlier successive linear programming framework of \citet{Ferris}, which also permits other trust-region norms. It corresponds formally to the case $\rho=1$ of the local linear optimization oracle defined by \citet{garberhazan2016}, whose construction on polytopes allows $\rho>1$. A cheap global linear oracle does not automatically make \eqref{eq:oracle} cheap. We assume exact access to this local oracle and count its calls separately from its implementation cost.

The radius determines how far the oracle can search from the current iterate. To obtain linear convergence, the radius must remain large enough relative to the current objective error to permit sufficient progress, while excessive radii must be reduced. Prescribing this relation directly typically requires growth or solution information that is unavailable to the algorithm. The original \algname{BroxGD} linear-rate radius, for example, depends on the distance to a solution \citep{locallmo}. We therefore ask:
\begin{quote}
Can we obtain a projection-free method for smooth strongly convex constrained minimization that converges linearly without knowing the optimal value, a solution, or the strong-convexity constant?
\end{quote}
\par\endgroup

\subsection{\texorpdfstring{Contributions}{Contributions}}
\begingroup\color{black}
We answer this question affirmatively by developing an adaptive short-step variant of \algname{BroxGD}. The analysis yields a more general result: linear convergence holds under quadratic growth, without requiring strong convexity or knowledge of the growth constant. Our method solves one local linear subproblem per update, minimizes a quadratic smoothness model along the resulting oracle segment, and chooses the next radius as the maximum of the accepted displacement and half the current radius. This adaptive rule maintains the required radius relation using quantities computed during each iteration. The method changes both the accepted point and the radius policy; it is a short-step variant of \algname{BroxGD}, not merely a new radius choice for its full-step update.

Our bound counts all local solves from any positive initial radius. It applies under quadratic growth, including strongly convex objectives and problems with multiple minimizers, and requires only a known smoothness bound $L$ as problem-dependent algorithmic input. The default $t_0=\|\nabla f(x_0)\|/L$ removes the initialization phase at a nonoptimal start (Corollary~\ref{cor:initialradius}). No projection onto $\X$ or function-value test is used by the outer method. Whether this oracle model yields a computational advantage depends on how the local subproblem is solved.
\par\endgroup

The central contribution is the combination of a \termrev{computable adaptive radius rule}, one exact local solve per update, and an explicit bound on all oracle calls needed for linear convergence under quadratic growth. The growth constant enters the analysis but is not supplied to the algorithm. The radius invariant explains how accepted movement supplies the scale that unavailable solution information would otherwise provide.
\begin{enumerate}
\item \textcolor{black}{\textbf{An adaptive method and an \termrev{oracle complexity bound}.} We establish linear convergence under quadratic growth on nonempty, closed, convex sets, including unbounded sets and problems with multiple minimizers. Each update uses one exact local linear-oracle call, with no projection onto $\X$, rejected trial solves, or function-value tests. The method requires $L$, but no growth constant, optimal value, minimizer, or solution distance. Theorem~\ref{thm:main} bounds all calls needed to reach objective error $\varepsilon$, explicitly tracking initialization and giving $O\!\left((L/\alpha)\log(1/\varepsilon)\right)$ for fixed problem and initialization.}
\item \textcolor{black}{\textbf{Radius control and computable initialization.} We prove an invariant linking the radius to the remaining objective error, and bound the time needed to establish it from any $t_0>0$ (Lemma~\ref{lem:envelope}). Once it holds, every update either achieves a specified objective contraction or halves the radius (Lemma~\ref{lem:dichotomy}). The default $t_0=\|\nabla f(x_0)\|/L$ establishes the invariant immediately at a nonoptimal start, eliminating the initialization phase without using solution information (Corollary~\ref{cor:initialradius}).}
\item \textcolor{black}{\textbf{Convergence of the entire trajectory.} Beyond the \termrev{iteration complexity bound}, we obtain an $R$-linear convergence bound for the objective values, finite length of the iterate sequence, and $R$-linear convergence of the iterates to a single solution (Corollary~\ref{cor:rlinear}). Appendix~\ref{app:lyapunov} gives a shorter Lyapunov proof after initialization, while explaining its less sharp constants and loss of the separate contraction and halving counts.}
\item \textcolor{black}{\textbf{Connections to standard methods.} A projection-dependent radius recovers \algname{ProjGD}, while a radius covering a compact feasible set recovers short-step \algname{FW} (Propositions~\ref{prop:pgd} and~\ref{prop:fw}). These are exact reductions of the local-oracle and segment step under alternative radius rules.}
\item \textcolor{black}{\textbf{What a fixed radius can and cannot guarantee.} For every fixed $t>0$, we prove an $O\!\left(Lt^2/k\right)$ objective-error bound and convergence in distance to the solution set (Theorem~\ref{thm:fixed}). Transferring a known \algname{FW} lower bound shows that, even for a fixed strongly convex problem, a constant radius need not admit a geometric bound uniform over all feasible initial points (Proposition~\ref{prop:lower}). This does not rule out geometric convergence for an individual initialization.}
\end{enumerate}
Under strong convexity, the bounds apply with the strong-convexity constant in place of the quadratic growth constant.

\subsection{Related work and oracle models}\label{sec:related}
\textcolor{black}{Table~\ref{tab:comparison} compares the number of oracle calls needed to reach objective error at most $\varepsilon$. The bounds show how this number grows as $\varepsilon$ decreases, for a fixed problem and initialization. Each row specifies the oracle being counted.}
\begin{table}[p]
\centering
\begin{threeparttable}\color{black}
\renewcommand{\rev}[1]{#1}
\hypersetup{linkcolor=black,citecolor=citationred}
\caption{\textcolor{black}{Total oracle calls to reach $f(x_k)-f_\star\le\varepsilon$, grouped by oracle model. A dash indicates that the cited result does not bound total calls; see the notes for accepted-step guarantees.}}
\label{tab:comparison}
\begingroup\fontsize{9}{10}\selectfont
\setlength{\tabcolsep}{3pt}
\renewcommand{\arraystretch}{1.0}
\newcommand{\tblnew}[1]{#1}
\definecolor{tableheadingblue}{HTML}{365F83}
\newcommand{\tblgroup}[1]{\multicolumn{4}{@{}l@{}}{\rule{0pt}{9pt}\textcolor{black}{\bfseries #1}}\\\specialrule{0.25pt}{1pt}{0pt}}
\newcommand{\tblmethod}[2]{\mbox{\textbf{#1}}\par\mbox{\citep{#2}}}
\newcommand{\tblinput}{\par\vspace{0pt}\textbf{Inputs:}\enspace}
\begin{tabular}{@{}>{\raggedright\arraybackslash}p{.27\textwidth}>{\raggedright\arraybackslash}p{.27\textwidth}>{\raggedright\arraybackslash}p{.20\textwidth}>{\raggedright\arraybackslash}p{\dimexpr.26\textwidth-18pt\relax}@{}}
\toprule
{\centering\textbf{Method}\par} & {\centering\textbf{Assumptions}\par\textit{Required inputs below}\par} & {\centering\textbf{Oracle / work}\par} & {\centering\textbf{Total calls}\par}\\
\midrule
\tblgroup{A. Projection and global linear oracles}
\tblmethod{\algnametiny{ProjGD}}{necoara2019}
& \tblnew{Convex; quadratic growth.}\tblinput $L$.
& One projection onto $\X$ per update.
& \tblnew{$O\!\left(\frac{L}{\alpha}\log\frac1\varepsilon\right)$\tnote{(a)}}\\
\specialrule{0.25pt}{1pt}{1pt}
\tblmethod{Restarted away-step \algnametiny{FW}}{kerdreux2022}
& \tblnew{Convex; quadratic growth.} \tblnew{\rev{Geometric scaling; finite away curvature.}}\tblinput \tblnew{Restart parameter.}
& \rev{Linear minimization over $\X$; active-set work and line search.}
& \tblnew{$O\left(\log\frac{1}{\varepsilon}\right)$\tnote{(b)}}\\
\midrule
\tblgroup{B. Exact local linear oracles}
\tblmethod{\tblnew{\algnametiny{SLP}}}{Ferris}
& \tblnew{\rev{Residual error bound; separated stationary values. Possibly nonconvex.}}\tblinput \tblnew{\rev{Trust-region parameters.}}
& \tblnew{\rev{Norm-ball local solves; trials and function-value tests.}}
& \tblnew{---\tnote{(c)}}\\
\specialrule{0.25pt}{1pt}{1pt}
\tblmethod{\algnametiny{BroxGD}}{locallmo}
& Strongly convex.\tblinput $L,\mu$, and $\|x_k-x_\star\|$.\tnote{(d)}
& Exact linear minimization over $\X\cap\B(x_k,t_k)$.
& \tblnew{$O\!\left(\frac{L}{\mu}\log\frac1\varepsilon\right)$\tnote{(d)}}\\
\specialrule{0.25pt}{1pt}{1pt}
\rowcolor{black!5}
\textbf{This work}
& \tblnew{Convex; quadratic growth.}\tblinput \rev{$L$; optional $t_0>0$. Default $t_0=\frac{\|\nabla f(x_0)\|}{L}$.}
& \rev{One exact solve over $\X\cap\B(x_k,t_k)$ per update. No function-value tests.}
& \tblnew{$O\!\left(\frac{L}{\alpha}\log\frac1\varepsilon\right)$\tnote{(e)}}\\
\midrule
\tblgroup{C. Expanded or approximate local oracles}
\tblmethod{Expanded local oracle}{garberhazan2016}
& Strongly convex. Polytope.\tblinput $L,\mu,\rho,H_0$.
& \tblnew{Linear minimization over $\X$ and vertex operations per call.}
& \tblnew{$O\!\left(\frac{L\rho^2}{\mu}\log\frac1\varepsilon\right)$\tnote{(f)}}\\
\specialrule{0.25pt}{1pt}{1pt}
\tblmethod{One-tangent method}{hazan2026}
& \tblnew{Strongly convex; compact body with interior rolling-ball radius $r_{\rm in}$.}\tblinput $L,\mu,G,H_0$ and geometric bounds.
& At most one tangent query per approximate local solve.
& $O\!\left(\left(\frac{L}{\mu}+\frac{G}{\mu r_{\rm in}}\right)\log\frac1\varepsilon\right)$\tnote{(g)}\\
\bottomrule
\end{tabular}
\begin{tablenotes}[flushleft]
\fontsize{8}{8.5}\selectfont\setlength{\itemsep}{0pt}\setlength{\parsep}{0pt}
\item[] \rev{\textbf{Convention.} Counts use column 3\textquotesingle s oracle. As $\varepsilon\downarrow0$, problem and initialization are fixed; \termrev{constants inside logarithms} and startup costs are suppressed.}
\item[(a)] \tblnew{Projection calls; conversion in Appendix~\ref{app:tablecurrency}.}
\item[(b)] \tblnew{The feasible set is compact. Linear oracle calls include drop steps; constants depend on geometry, support, and initialization.}
\item[(c)] \tblnew{\algnametiny{SLP}: successive linear programming. The cited rate counts accepted iterations, excludes trial solves, and allows a nonglobal limit; no stated-target call bound.}
\item[(d)] \rev{The first version's radius depends on $\|x_k-x_\star\|$; the displayed rate counts exact local calls for that rule (conversion in Appendix~\ref{app:tablecurrency}). A later revision also added radius backtracking \citep{broxgd2026v2}. Our results were established by early September 2026, before that development.}
\item[(e)] \tblnew{Exact local calls; see \eqref{eq:nonasymptotic} for initialization and radius dependence.}
\item[(f)] \tblnew{Expanded local oracle calls.}
\item[(g)] \tblnew{Exact supporting-tangent queries, not membership queries.}
\item[] \textbf{Notation.} $f$: smooth objective; $\X$: nonempty closed convex feasible set; $f_\star=\min_\X f$; $\Xs=\arg\min_\X f$; $x_\star\in\Xs$: minimizer; $x_k$: iterate $k$ ($x_0$: initial point); $\varepsilon>0$: objective-error tolerance; $\nabla f$: gradient; $\|\cdot\|$: Euclidean norm; $\B(x,t)=\left\{y:\|y-x\|\le t\right\}$; $t_k$: radius at iteration $k$ ($t_0$: initial radius).
\item[] \textbf{Curvature.} $L>0$: smoothness bound (\eqref{eq:smooth} for this work); $\mu>0$: strong-convexity constant. For \algnametiny{ProjGD} and this work, $\alpha=\min\left\{\alpha_0,L\right\}$, where $f(x)-f_\star\ge\frac{\alpha_0}{2}\operatorname{dist}(x,\Xs)^2$ on $\X$; $\alpha_0>0$ is the quadratic-growth constant and $\operatorname{dist}$ is Euclidean distance. Neither method needs $\alpha_0$ as input.
\item[] \textbf{Other symbols.} $H_0\ge f(x_0)-f_\star$: initial error bound; $\rho\ge1$: oracle radius-expansion factor; $G\ge\sup_{x\in\X}\|\nabla f(x)\|$; $r_{\rm in}>0$: uniform interior rolling-ball radius; $\log$: natural logarithm; $O$: asymptotic upper order. \algnametiny{ProjGD}: projected gradient descent; \algnametiny{FW}: Frank--Wolfe.
\item[] \rev{\textbf{Cost caveat.} Ancillary operations are listed but not counted. Different oracle types need not have equal computational costs.}
\end{tablenotes}
\endgroup
\end{threeparttable}
\end{table}

\begingroup\color{black}
\paragraph{Ball optimization oracles.}
Local minimization over balls predates the Broximal framework. \citet{balloracle2020} study an oracle that minimizes the convex objective itself over a Euclidean ball and develop an accelerated outer method with matching oracle lower bounds up to logarithmic factors. \citet{insideball2021} refine this approach and implement ball oracles for softmax objectives arising in maximal-loss minimization. \citet{ballgame2023} extend ball-oracle acceleration to non-Euclidean geometry and apply it to matrix games and maxima of smooth functions. These works separate the cost of implementing a ball subproblem from the number of outer oracle calls. Their objective-minimization oracle differs from \eqref{eq:oracle}, which minimizes a linear model over $\X\cap\B(x_k,t_k)$; their complexity guarantees therefore do not directly give the local-linear-oracle bound studied here.

\paragraph{The Broximal framework.}
The ball-proximal, or Broximal, point method of \citet{bpm2025} replaces the proximal quadratic penalty by a ball constraint and studies repeated minimization of the objective over that ball. Extensions include arbitrary norm geometries \citep{nonEuclideanBrox2025} and global nonconvex convergence under Broximal Alignment \citep{alignment2026}. The trust-region stabilized proximal method of \citet{trppm2026} connects proximal regularization and ball constraints, with an equivalence to the Broximal point method in its active-constraint regime. \algname{BroxGD} instead uses a linear model in the local subproblem. Our contribution is an \termrev{implementable adaptive radius rule} combined with a short step along its oracle segment.

Several further developments appeared after our results were established in early September 2026. \citet{sharpbrox2026} sharpen the theory of exact ball-proximal minimization and analyze radius sequences and adaptive radius choices. Their subproblem minimizes the original objective, whereas our rule is analyzed for linear-model subproblems. \citet{projectionbrox2026} give a projection interpretation of local linear minimization, and \citet{fwbrox2026} interpret Frank--Wolfe as an approximation of \algname{BroxGD}; we relate these perspectives to our reductions in Section~\ref{sec:radius}. \citet{inexactbrox2026} analyze perturbations of ball-proximal points, providing complementary results for inexact objective-minimization subproblems.
\endgroup

\paragraph{Quadratic growth and projected methods.}
\citet{drusvyatskiylewis2018} relate quadratic growth on sublevel sets to residual error bounds for proximal methods. \citet[Section~5.1.2]{necoara2019} establish linear convergence of \algname{ProjGD} under quadratic growth. Appendix~\ref{app:pgdqg} derives the objective contraction used in our comparison.

\paragraph{Successive linear programming.}
\citet{Ferris} minimize linear functions over a closed convex set intersected with a norm ball. In the Euclidean case, their subproblem is \eqref{eq:oracle}. They also relate its solutions to suitable projected-gradient displacements. Their successive linear programming algorithm enlarges or shrinks the trial radius using function-decrease tests. Building on \citet{luotseng1993}, they prove asymptotic $Q$-linear convergence of the objective values under a residual error bound and a condition separating stationary values, including for nonconvex objectives. Our method uses one exact local solve per update without function-value tests, and our bound counts every solve. The successive linear-quadratic programming method of \citet[Algorithm~2.1, Step~5]{byrd2005} also uses movement-based radius control. It relates the next linear-subproblem radius to the Cauchy displacement within a scheme that uses two trust regions and acceptance tests. In one branch of \algname{NEWUOA}, \citet{powell2004} sets the radius to the maximum of the step norm and half the previous radius. The ratio of actual to predicted decrease selects this branch, and a separate radius safeguard follows.

\begingroup\color{black}
\paragraph{Exact local oracles: full steps and short steps.}
The solution-dependent radius in the original \algname{BroxGD} paper \citep[Theorem~3.5]{locallmo} gives a sharp distance contraction under strong convexity; Appendix~\ref{app:tablecurrency} converts it to objective accuracy. We retain that version-specific benchmark in Table~\ref{tab:comparison}. Our adaptive short-step method and its convergence results were developed independently of the adaptive variants added in the September 2026 revision of that paper. All our results were established by early September 2026, before those developments, and inspired the later adaptive variants. In that later revision, \citet[Algorithm~2 and Theorem~3.6]{broxgd2026v2} instead backtrack from $\|\nabla f(x_k)\|/L$, shrinking the radius by $\xi\in(0,1)$ until the oracle displacement $d$ satisfies $-\inner{\nabla f(x_k)}{d}\ge L\|d\|^2$. The accepted point is the full oracle point. With a globally Lipschitz gradient and strong convexity, the objective error is bounded by $\left(1+\xi^2\mu/L\right)^{-K}$ times its initial value. Here $K$ counts accepted steps; rejected trials cost additional calls. Neither the test nor its implementation needs $\mu$ or function values. That revision also proves a geometric-radius guarantee from certified initial bounds.

Our clipped step avoids a search for an acceptable full oracle point. \termrev{Its radius update uses} the accepted displacement, and our bound counts all solves under quadratic growth and the \termrev{quadratic upper-bound assumption} \eqref{eq:smooth}. These distinctions concern the algorithm, assumptions, and counted resource; they do not establish runtime superiority. The projected-gradient condition discussed in Remark~\ref{rem:projectedpl} refers specifically to the original version.
\endgroup

\paragraph{Expanded local oracles.}
\citet{garberhazan2016} construct a local linear optimization oracle on polytopes using one global linear minimization and operations on a vertex decomposition. The returned point may be as far as $\rho t$ from the current iterate, where $\rho\ge1$, but its linear objective value is no larger than the minimum over $\X\cap\B(x,t)$. For strongly convex objectives, their algorithm chooses shrinking radii using smoothness, strong convexity, $\rho$, and an initial error bound. Our results instead assume direct access to exact minimization over $\X\cap\B(x,t)$ for nonempty, closed, convex sets. This corresponds formally to $\rho=1$. \citet{dvurechensky2020} use the expanded oracle for self-concordant objectives on polytopes, with step-dependent radius shrinkage and a positive curvature lower bound on the relevant level set.

\paragraph{Restarted FW and geometric scaling.}
\citet{kerdreux2022} combine a H\"olderian error bound with a positive geometric scaling inequality. Under quadratic growth, their method requires a logarithmic total number of global linear oracle calls, including drop steps. The bound depends on the scaling constant and away curvature. Pyramidal width supplies a positive scaling constant on polytopes. Their results beyond polytopes require additional conditions, such as an interior solution or a strongly convex feasible set, depending on the regime.

\paragraph{Other adaptive mechanisms.}
Away-step methods achieve linear convergence on polytopes for smooth objectives $f(x)=h(Ax)+b^\top x$ with strongly convex $h$, even when $f$ is not strongly convex \citep{lacostejulien2015,beckshtern2017}. \citet[Theorems~5 and~11]{garberwolf2021} obtain linear convergence under quadratic growth on polytopes with a nearest-extreme-point oracle and correction steps. Their adaptive parameter search incurs additional oracle work. \citet{rinaldizeffiro2023} chain short steps while keeping the gradient fixed. Their analysis assumes a directional angle condition and a constrained Kurdyka--\L{}ojasiewicz inequality with exponent $\frac{1}{2}$. For away-step and pairwise \algname{FW}, their Proposition~4.3 bounds the worst-case average number of inner steps by two along an infinite sequence, when the linear minimizer is fixed within each chain. \citet{pedregosa2020} adapt smoothness estimates through backtracking and allow approximate global linear subproblems. \citet[Theorems~4.2 and~4.4]{takahashi2026} use adaptive Bregman steps; their convex linear-rate results require an interior minimizer or polytope geometry, together with the stated error-bound and Bregman regularity assumptions.

\paragraph{Geometric implementations.}
\citet{hazan2026} implements an approximate local oracle through supporting-tangent or membership queries for smooth strongly convex objectives on compact convex bodies satisfying a uniform interior rolling-ball condition. The guarantee depends on geometric and gradient bounds. Our analysis assumes an exact local oracle, with $L$ as its only required problem-dependent algorithmic parameter.

\paragraph{Lower bounds.}
\citet{halbey2026} give a fixed strongly convex, ball-constrained problem with horizon-dependent initializations that obstruct a uniform geometric bound for \algname{FW} with short steps or exact line search. \citet{grimmerliu2026} give a lower bound for a broader deterministic first-order and global linear oracle model using horizon-dependent hard sets. Proposition~\ref{prop:lower} transfers the \algname{FW} obstruction to our fixed-radius variant when the ball covers $\X$. Neither result gives a lower bound for our method.

\Needspace{5\baselineskip}
\section{Algorithm and assumptions}\label{sec:algorithm}
\begingroup\color{black}
We specify the problem and oracle, derive the accepted step from the smoothness model, and then explain how that step determines the next radius.
\endgroup

\subsection{\texorpdfstring{Problem and oracle}{Problem and oracle}}\label{sec:assumptions}
Throughout the analysis, norms and projections are Euclidean. We assume that $f$ is differentiable on an open set containing $\X$, convex on $\X$, and satisfies
\begin{equation}\label{eq:smooth}
 f(y)\le f(x)+\inner{\nabla f(x)}{y-x}+\frac L2\|y-x\|^2
 \qquad(x,y\in\X)
\end{equation}
for a known $L>0$. We assume that the minimum is attained and write
\[
 f_\star=\min_{x\in\X}f(x),\qquad
 \Xs=\left\{x\in\X:f(x)=f_\star\right\}.
\]
The set $\Xs$ is nonempty, closed, and convex. Thus the distance $\dist(x,\Xs)=\inf_{y\in\Xs}\|x-y\|$ is attained at a nearest solution. We write $x_\star\in\Xs$ for a minimizer, without assuming uniqueness.

\begin{assumption}[Quadratic growth]\label{ass:qg}
There is $\alpha_0>0$ such that
\begin{equation}\label{eq:qg}
 f(x)-f_\star\ge\frac{\alpha_0}{2}\dist(x,\Xs)^2
 \qquad\text{for every }x\in\X.
\end{equation}
\end{assumption}
\begin{remark}[Weaker sufficient domains]\label{rem:qgdomain}
The objective values do not increase, so quadratic growth on the initial sublevel set $\left\{x\in\X:f(x)\le f(x_0)\right\}$ also suffices. In fact, the proofs use quadratic growth only at the iterates. To cover every choice of exact oracle minimizer, the same growth constant must apply to all resulting trajectories. Assumption~\ref{ass:qg} ensures this for each radius variant considered below.
\end{remark}
For the analysis, set
\begin{equation}\label{eq:scales}
 \alpha=\min\left\{\alpha_0,L\right\}.
\end{equation}
Then $\dist(x_k,\Xs)\le\sqrt{\frac{2\left(f(x_{k})-f_\star\right)}{\alpha}}$ and $0<\alpha\le L$. The algorithm does not evaluate $\alpha$, $f_\star$, or this distance bound.

\begingroup\color{black}
For each $x_k\in\X$ and $t_k>0$, evaluate $g_k=\nabla f(x_k)$, choose any exact minimizer $z_k$ in \eqref{eq:oracle}, and define
\begin{equation}\label{eq:stepnotation}
 d_k=z_k-x_k,\qquad D_k=\|d_k\|,\qquad
 G_k=-\inner{g_k}{d_k}.
\end{equation}
Comparison with the feasible point $x_k$ gives $G_k\ge0$. In fact, $G_k=0$ is equivalent to global optimality (Lemma~\ref{lem:descent}), so the algorithm stops in that case. Otherwise $G_k>0$ and $D_k>0$, and all divisions below are well defined. Exactness is used both for this certificate and for the comparison estimates in the proof.

\subsection{\texorpdfstring{Choosing the step}{Choosing the step}}\label{sec:adaptive-method}
The oracle gives a feasible direction. To choose how far to move, write a segment point as $x_k+\eta d_k$, with $\eta\in[0,1]$. The upper model is
\begin{equation}\label{eq:segmentmodel}
 f(x_k+\eta d_k)\overset{\eqref{eq:smooth},\eqref{eq:stepnotation}}{\le}
 f(x_k)-\eta G_k+\frac{L\eta^2D_k^2}{2}.
\end{equation}
Minimizing the right-hand side over $[0,1]$ gives the clipped coefficient and accepted displacement
\begin{equation}\label{eq:coefficient}
 \eta_k=\min\left\{1,\frac{G_k}{LD_k^2}\right\},
 \qquad x_{k+1}=x_k+\eta_k d_k.
\end{equation}
Its length is therefore
\begin{equation}\label{eq:length}
 \ell_k:=\|x_{k+1}-x_k\|
 \overset{\eqref{eq:coefficient},\eqref{eq:stepnotation}}{=}
 \eta_kD_k
 \overset{\eqref{eq:coefficient}}{=}
 \min\left\{\frac{G_k}{LD_k},D_k\right\}.
\end{equation}
Let $S_k=[x_k,z_k]$. Completing the square gives an equivalent segment-projection interpretation:
\begin{equation}\label{eq:segment}
\begin{relchain}
x_{k+1} &\overset{\eqref{eq:coefficient}}{=}&
 \argmin_{v\in S_k}
 \left\{f(x_k)+\inner{g_k}{v-x_k}+\frac L2\|v-x_k\|^2\right\}\\[4pt]
 & = & \Proj_{S_k}\left(x_k-\frac{g_k}{L}\right).
\end{relchain}
\end{equation}
Once $z_k$ is available, only scalar clipping is needed. The method is projection-free at the outer level: it never projects onto $\X$. This does not impose an implementation strategy on the local-oracle solver.

\subsection{\texorpdfstring{Updating the radius}{Updating the radius}}\label{sec:radius-update}
The accepted length supplies a scale already supported by model descent. We retain that scale for the next local search, but do not let a single small displacement shrink the radius by more than a factor of two:
\begin{equation}\label{eq:update}
 x_{k+1}\overset{\eqref{eq:coefficient},\eqref{eq:length}}{=}x_k+\ell_k\frac{d_k}{D_k},
 \qquad t_{k+1}=\max\left\{\ell_k,\frac{t_k}{2}\right\}.
\end{equation}
Since $\ell_k\le D_k\le t_k$, the radius stays in $[t_k/2,t_k]$. A full step that reaches the ball boundary keeps the radius unchanged. A short accepted displacement permits reduction, while the safeguard prevents an abrupt collapse. Section~\ref{sec:convergence} shows how these choices maintain a lower bound on the radius relative to the remaining error, although the algorithm cannot evaluate that error.

We distinguish three step types. A \emph{full boundary step} has $\ell_k=D_k=t_k$. A \emph{full interior step} has $\ell_k=D_k<t_k$. A \emph{short step} has $\ell_k<D_k$. \termrev{Here boundary and interior refer to the local ball, not to the feasible set $\X$.} Equality in the clipping test belongs to the full-step case. Algorithm~\ref{alg:main} collects the updates; its supplied-radius and default-radius branches generate the same method after initialization.

\Needspace{15\baselineskip}
\begingroup\color{black}
\begin{example}[\termrev{Two updates illustrating the radius update}]\label{ex:controller}
Take $f(x)=x^2/2$ on $\X=[0,2]$, the \termrev{valid quadratic upper-bound constant} $L=2$, and $(x_0,t_0)=(1,1/2)$. The increasing local linear objective selects the left endpoint of the local interval. The first two oracle points are $z_0=1/2$ and $z_1=0$, with $D_0=D_1=1/2$, $G_0=1/2$, and $G_1=1/4$. Hence
\[
\begin{aligned}
 \ell_0&\overset{\eqref{eq:length}}{=}\min\left\{\frac12,\frac12\right\}=\frac12,
 & (x_1,t_1)&\overset{\eqref{eq:update}}{=}\left(\frac12,\frac12\right),\\[3pt]
 \ell_1&\overset{\eqref{eq:length}}{=}\min\left\{\frac14,\frac12\right\}=\frac14,
 & (x_2,t_2)&\overset{\eqref{eq:update}}{=}\left(\frac14,\frac14\right).
\end{aligned}
\]
The full boundary step preserves the radius; the next short step halves it. Both updates contract the objective error, so halving need not indicate failed contraction. Section~\ref{sec:convergence} proves the general guarantee. Figure~\ref{fig:controller} illustrates both updates.
\end{example}
\par\endgroup

\begin{figure}[H]
\centering
\includegraphics[width=\textwidth]{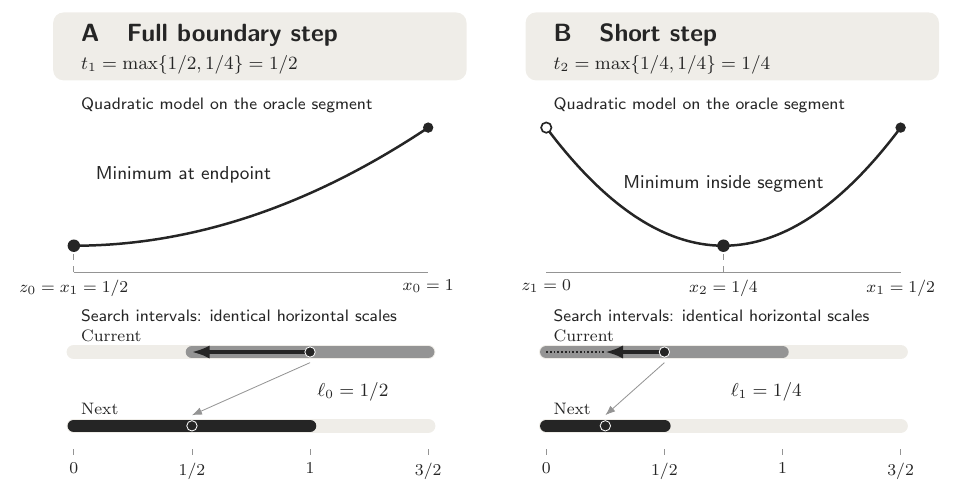}
\begingroup\color{black}
\caption{Two consecutive updates in Example~\ref{ex:controller}. Top: minimizing the quadratic upper model along the oracle segment accepts its endpoint in A and an interior point in B; the vertical scales differ. Bottom: the current local interval (light gray) and next local interval (dark gray) share a common horizontal scale. Dots mark their centers; horizontal arrows show accepted movement. The first accepted length equals the radius, so the radius stays at $1/2$. The second length is $1/4$, so the radius halves. Both updates reduce the objective error by a factor of four.}
\label{fig:controller}
\endgroup
\end{figure}

\begin{algorithm}[H]
\color{black}
\caption{\textcolor{black}{Adaptive Short-Step Broximal Gradient Descent}}\label{alg:main}
\begin{algorithmic}[1]
\Require $x_0\in\X$, $L>0$ satisfying \eqref{eq:smooth}, and optionally $t_0>0$.
\If{no initial radius is supplied}
  \State Evaluate $g_0=\nabla f(x_0)$.
  \If{$g_0=0$}
    \State \Return $x_0$ and terminate.
  \EndIf
  \State Set $t_0=\|g_0\|/L$.
\EndIf
\For{$k=0,1,\ldots$}
  \State Evaluate $g_k=\nabla f(x_k)$ if not already available.
  \State Choose any
  \[
    z_k\in\argmin_{z\in\X\cap\B(x_k,t_k)}\inner{g_k}{z}.
  \]
  \State Set $d_k=z_k-x_k$, $D_k=\|d_k\|$, and $G_k=-\inner{g_k}{d_k}$.
  \If{$G_k=0$}
    \State \Return $x_k$ and terminate.
  \EndIf
  \State Compute $\ell_k=\min\left\{\frac{G_k}{LD_k},D_k\right\}$.
  \State Update the iterate and radius:
  \[
    x_{k+1}=x_k+\ell_k\frac{d_k}{D_k},
    \qquad t_{k+1}=\max\left\{\ell_k,\frac{t_k}{2}\right\}.
  \]
\EndFor
\end{algorithmic}
\end{algorithm}
Each update uses one gradient evaluation and one exact local-oracle call. There are no rejected trial solves and no function-value evaluations. The prescribed objective accuracy is an analysis target; the zero-gap test certifies exact optimality, not a computable positive-error tolerance.
\endgroup

\Needspace{10\baselineskip}
\section{Convergence under quadratic growth}\label{sec:convergence}
\begingroup\color{black}
The proof tracks the objective gap and \textcolor{black}{an upper bound on the distance to the solution set}, neither of which the algorithm evaluates:
\begin{equation}\label{eq:analysisnotation}
 \Delta_k=f(x_k)-f_\star,\qquad s_k=\sqrt{\frac{2\Delta_k}{\alpha}}.
\end{equation}
Quadratic growth gives $\dist(x_k,\Xs)\le s_k$. The radius invariant is $t_k\ge\alpha s_k/(2L)$. Its role is to keep the search radius large enough relative to this upper bound on solution distance.

The mechanism has three stages. Before the invariant holds, every step is full: boundary steps keep the radius fixed while decreasing $s_k$, and an interior step establishes the invariant immediately. Once it holds, every update either contracts $\Delta_k$ by the factor $1-\alpha/(4L)$ or halves an excessive radius. Counting initialization steps, contractions, and the remaining halvings gives a bound on every local-oracle call. The analysis does not claim that each update contracts the objective error by the same factor.

We state the bound first, then establish the descent and radius-control lemmas. Appendix~\ref{app:count} gives the complete count, and Appendix~\ref{app:fullsequence} derives \textcolor{black}{an $R$-linear convergence bound for the objective values} and convergence to a single solution. The shorter potential proof in Appendix~\ref{app:lyapunov} provides an alternative after initialization, with less sharp constants. Write $[v]_+=\max\left\{v,0\right\}$. Throughout, $\log$ is the natural logarithm and $\log_2$ is the base-two logarithm.
\endgroup

\begin{theorem}[Total exact-oracle complexity]\label{thm:main}
Under the assumptions of Section~\ref{sec:algorithm}, let $x_0\in\X$ and $t_0>0$. If $x_0$ is optimal, the algorithm terminates after at most one local oracle call. Otherwise, for $0<\varepsilon<\Delta_0$, let $K(\varepsilon)$ be the first index satisfying $\Delta_{K(\varepsilon)}\le\varepsilon$. Then
\begin{equation}\label{eq:nonasymptotic}
 K(\varepsilon)\le
 \frac4{t_0}\pos{s_0-\frac{2Lt_0}{\alpha}}
 +\frac{4L}{\alpha}\log\frac{\Delta_0}{\varepsilon}
 +\pos{\log_2\frac{t_0\sqrt{\alpha}}{\sqrt{2\varepsilon}}}+3.
\end{equation}
Consequently, for fixed initialization and problem constants, we obtain
\begin{equation}\label{eq:asymptotic}
 K(\varepsilon)=O\!\left(\frac L\alpha\log\frac1\varepsilon\right)
 \qquad(\varepsilon\downarrow0).
\end{equation}
Each update uses one exact local oracle call. If an optimum is attained, checking the stopping condition there can require one additional call.
\end{theorem}
Appendix~\ref{app:count} gives the sharper bound \eqref{eq:complexity}, which counts the three types of update separately. The first term of \eqref{eq:nonasymptotic} quantifies the initialization phase for a small $t_0$. It vanishes as soon as $t_0\ge\frac{\alpha}{2L}s_0$. Enlarging $t_0$ beyond this threshold affects only the logarithmic halving term in the bound. An excessively small radius can therefore lengthen the initialization phase, whereas a large radius enters the bound only logarithmically. The following choice meets the threshold without requiring its value.

\begin{corollary}[Gradient-based initial radius]\label{cor:initialradius}
Under the assumptions of Theorem~\ref{thm:main}, if $\nabla f(x_0)=0$, then $x_0$ is optimal. Otherwise choose
\begin{equation}\label{eq:initialradius}
 t_0=\frac{\|\nabla f(x_0)\|}{L}.
\end{equation}
Whenever the starting point is nonoptimal, this choice establishes the radius invariant \eqref{eq:invariant} immediately, so no initialization phase is needed and
\begin{equation}\label{eq:initialcomplexity}
 K(\varepsilon)\le
 \simplerev{\left\lceil\frac{4L}{\alpha}\log\!\left(\frac{\Delta_0}{\varepsilon}\right)\right\rceil}
 +\ceilp{\log_2\frac{\|\nabla f(x_0)\|\sqrt{\alpha}}{L\sqrt{2\varepsilon}}}.
\end{equation}
The initial gradient is reused, so no additional local oracle call is needed.
\end{corollary}
\begin{proof}[\normalfont\bfseries Proof]
At a nonoptimal $x_0$, convexity, Cauchy-Schwarz and quadratic growth imply
\[
 0<\Delta_0\le\|\nabla f(x_0)\|\dist(x_0,\Xs)
 \overset{\textcolor{black}{\eqref{eq:qg},\eqref{eq:scales}}}{\textcolor{black}{\le}}\|\nabla f(x_0)\|s_0.
\]
Dividing by $s_0>0$ gives
\[
 \|\nabla f(x_0)\|\overset{\eqref{eq:analysisnotation}}{\ge}\frac{\alpha}{2}s_0.
\]
Thus $t_0=\frac{\|\nabla f(x_0)\|}{L}$ satisfies the radius invariant at the starting point. The initialization term in the sharper count \eqref{eq:complexity} is zero, and substituting this radius, \simplerev{then using $1/[-\log(1-\alpha/(4L))]\le4L/\alpha$,} gives the stated bound. If the initial gradient is zero, convexity gives optimality directly.
\end{proof}
The local zero-gap test also certifies optimality at constrained solutions with a nonzero gradient. A small positive local gap alone does not certify $\Delta_k\le\varepsilon$ without additional bounds relating the gap to the objective error. Thus $K(\varepsilon)$ is a convergence bound, not a stopping rule that the algorithm can evaluate from its inputs.

\subsection{Descent and local comparison}
The following lemma establishes that the method is well defined, characterizes optimality through the local gap, and bounds the objective decrease. The zero-gap characterization is also given by \citet{Ferris}. We include a proof in the appendix for completeness.
\begin{lemma}[Well-posedness and model decrease]\label{lem:descent}
The local oracle is well defined, and $\inner{\nabla f(x_{k})}{x_{k}-z_{k}}=0$ if and only if $x_k$ is optimal. At each nonterminal iteration,
\begin{equation}\label{eq:feasible}
 x_{k+1}\in\X,\qquad 0<\ell_k\le\|z_{k}-x_{k}\|\le t_k,
 \qquad \frac{t_k}{2}\le t_{k+1}\le t_k.
\end{equation}
Moreover,
\begin{equation}\label{eq:descent}
 \Delta_k-\Delta_{k+1}\ge
 \ell_k\frac{\inner{\nabla f(x_{k})}{x_{k}-z_{k}}}{\|z_{k}-x_{k}\|}-\frac L2\ell_k^2
 \ge\frac{\ell_k}{2}\frac{\inner{\nabla f(x_{k})}{x_{k}-z_{k}}}{\|z_{k}-x_{k}\|}>0.
\end{equation}
A short step gives
\begin{equation}\label{eq:short}
 \Delta_k-\Delta_{k+1}\ge\frac{\inner{\nabla f(x_{k})}{x_{k}-z_{k}}^2}{2L\|z_{k}-x_{k}\|^2},
\end{equation}
and a full step gives
\begin{equation}\label{eq:full}
 \Delta_k-\Delta_{k+1}\ge\inner{\nabla f(x_{k})}{x_{k}-z_{k}}-\frac L2\|z_{k}-x_{k}\|^2
 \ge\frac{1}{2}\inner{\nabla f(x_{k})}{x_{k}-z_{k}}.
\end{equation}
\end{lemma}
To compare the oracle point with any feasible $y$, we move from $x_k$ toward $y$ only as far as the ball $\B(x_k,t_k)$ allows. The next lemma quantifies this comparison.
\begin{lemma}[Scaling]\label{lem:scaling}
For $y\in\X\setminus\left\{x_k\right\}$,
\begin{equation}\label{eq:scaling}
 \inner{\nabla f(x_{k})}{x_{k}-z_{k}}\ge\min\!\left\{1,\frac{t_k}{\|y-x_k\|}\right\}\inner{\nabla f(x_k)}{x_k-y}.
\end{equation}
At a nonoptimal iterate, comparison with a nearest solution yields
\begin{equation}\label{eq:slopes}
\begin{relchain}
\inner{\nabla f(x_{k})}{x_{k}-z_{k}} & \ge & \min\!\left\{1,\frac{t_k}{\dist(x_k,\Xs)}\right\}\Delta_{k},\\[4pt]
\frac{\inner{\nabla f(x_{k})}{x_{k}-z_{k}}}{\|z_{k}-x_{k}\|} & \ge & \frac{\Delta_k}{\max\left\{t_k,\dist(x_k,\Xs)\right\}}.
\end{relchain}
\end{equation}
\end{lemma}
The distance $\dist(x_k,\Xs)$ from the current iterate to the solution set is bounded above by $s_k$ under quadratic growth, without requiring a diameter bound.

If the oracle solution lies strictly inside the ball, it also minimizes the linear objective over the whole feasible set. This property gives the following bound on the next radius.
\begin{lemma}[\termrev{Objective-error bound after a full interior step}]\label{lem:reset}
If a full interior step occurs, then $z_k$ minimizes $\inner{\nabla f(x_k)}{z}$ over all of $\X$, and
\begin{equation}\label{eq:reset}
 \Delta_{k+1}\le\frac L2\|z_{k}-x_{k}\|^2,
 \qquad t_{k+1}\ge\|z_{k}-x_{k}\|\ge\sqrt{\frac{2\Delta_{k+1}}L}.
\end{equation}
\end{lemma}
Because $\|z_k-x_k\|<t_k$, a sufficiently short move from $z_k$ toward any $y\in\X$ remains inside the ball and, by convexity, inside $\X$. If $y$ had a smaller linear objective value than $z_k$, this move would contradict local oracle optimality. Hence $z_k$ also minimizes the linear objective on all of $\X$.

\subsection{Radius control}
\begingroup\color{black}
We call $t_k\ge\alpha s_k/(2L)$ an invariant because it is preserved by every subsequent update once established. The three relevant regimes are $t_k<\alpha s_k/(2L)$ (initialization), $\alpha s_k/(2L)\le t_k\le s_k$ (guaranteed contraction), and $t_k>s_k$ (contraction or radius reduction). The next two lemmas justify this description.
\endgroup

\begingroup\color{black}
\begin{lemma}[Invariant and initialization]\label{lem:envelope}
The radius bound
\begin{equation}\label{eq:invariant}
 t_k\ge\frac{\alpha}{2L}s_k
\end{equation}
is preserved from one iterate to the next once it is satisfied. Starting from any $t_0>0$, \eqref{eq:invariant} holds for some index $k=N$ satisfying
\begin{equation}\label{eq:catchup}
 N\le C_0:=\ceilp{\frac4{t_0}\left(s_0-\frac{2Lt_0}{\alpha}\right)}.
\end{equation}
If $t_k<\frac{\alpha}{2L}s_k$, a short step is impossible. A full interior step establishes the invariant at the next iterate. Otherwise, a full boundary step satisfies
\begin{equation}\label{eq:catchupstep}
 t_{k+1}=t_k,\qquad s_k-s_{k+1}\ge \frac{t_k}{4}.
\end{equation}
\end{lemma}
Appendix~\ref{app:invariant} shows how each of the three step types preserves this bound. Before the bound holds, every step is full. Boundary steps keep the radius at $t_0$ and reduce the upper bound on the distance to the solution set by at least $\frac{t_0}{4}$ per update. An interior step establishes the invariant at the next iterate and ends the initialization phase.

\begin{lemma}[Objective contraction or radius halving]\label{lem:dichotomy}
Suppose $t_k\ge\frac{\alpha}{2L}s_k$ at a nonterminal iteration. Either the objective error contracts according to
\begin{equation}\label{eq:contract}
 \Delta_{k+1}\le\left(1-\frac{\alpha}{4L}\right)\Delta_k,
\end{equation}
or the step is short, $t_k>s_k$, $\ell_k\le \frac{t_k}{2}$, and the radius update gives
\begin{equation}\label{eq:halve}
 t_{k+1}=\frac{t_k}{2}.
\end{equation}
\end{lemma}
These two options account for every update once the invariant holds. If a step does not achieve the stated contraction, it still decreases $f$ and halves the radius. Such a step can occur only when $t_k>s_k$. This stronger radius condition limits the number of such halvings before the target objective accuracy is reached. This gives the halving term in the complexity bound.
\endgroup

\section{The role of the radius}\label{sec:radius}
A fixed radius is a natural alternative to the adaptive rule. This raises the question of whether radius adaptation is needed for a general linear-convergence guarantee. To address it, we study the convergence of the constant-radius variant and its limitations. We first relate the local segment step to two familiar methods. A projection-dependent radius gives \textcolor{black}{\algname{ProjGD}}, while a radius covering a compact feasible set gives short-step \algname{FW}. Throughout this section, we vary the radius rule while keeping the local oracle and segment step unchanged.
Subsequent work by \citet{projectionbrox2026} identifies the full local-oracle point, for radii no larger than the Polyak radius, as a projection of the current iterate onto the feasible set intersected with a separating half-space. This differs from our segment projection \eqref{eq:segment}, which chooses the accepted point after the local oracle has been solved. \citet{fwbrox2026} obtain Frank--Wolfe by replacing the ball intersection with a smaller homothetic copy of the feasible set. Our Frank--Wolfe reduction below instead keeps the oracle exact and chooses a radius covering the whole feasible set.

\subsection{Two exact reductions}
{We first recover short-step \algname{FW} by choosing a radius that covers the feasible set, then recover \algname{ProjGD} through a projection-dependent radius.}

\begin{proposition}[Short-step Frank--Wolfe]\label{prop:fw}
If $\X$ is compact with diameter $D$ and $t_k=t\ge D$ for all $k$, each resulting update is a short-step \algname{FW} update.
\end{proposition}
The local feasible set then equals $\X$. Writing $x_{k+1}=x_k+\eta_k(z_k-x_k)$ gives
\[
 \eta_k\overset{\eqref{eq:coefficient},\eqref{eq:stepnotation}}{=}\min\!\left\{\frac{\inner{\nabla f(x_{k})}{x_{k}-z_{k}}}{L\|z_{k}-x_{k}\|^2},1\right\},
\]
the usual clipped \algname{FW} short-step coefficient \citep{pedregosa2020}. This is the coefficient in their formulation with smoothness estimate $L$ and maximum segment coefficient one.

The connection between Euclidean projection and local linear minimization is classical \citep{Ferris}. \citet{locallmo} use the projection-dependent radius below. The following proposition shows that the projected point is the unique local oracle solution at this radius and that our segment step reaches it.
\begin{proposition}[Projected gradient descent]\label{prop:pgd}
Fix $0<\gamma\le\frac{1}{L}$ and let
\begin{equation}\label{eq:pgdradius}
 p_k=\Proj_\X(x_k-\gamma\nabla f(x_k)),\qquad t_k=\|p_k-x_k\|.
\end{equation}
If $t_k=0$, then $x_k$ is optimal. Otherwise $p_k$ is the unique local oracle solution at this radius, and the segment step accepts $p_k$. Applying \eqref{eq:pgdradius} at each iteration therefore gives \textcolor{black}{\algname{ProjGD}} with stepsize $\gamma$.
\end{proposition}
The radius in \eqref{eq:pgdradius} is the norm of the projected-gradient displacement. Evaluating this radius directly requires the projection itself. The identity therefore explains the connection to \textcolor{black}{\algname{ProjGD}}, while Algorithm~\ref{alg:main} computes its radii from accepted step lengths and previous radii.

\subsection{Fixed-radius convergence and its limitation}
{With a constant radius, the method still converges, but its general guarantee is sublinear and a uniform geometric bound need not hold.}

\begin{theorem}[Fixed-radius upper bound]\label{thm:fixed}
Under the assumptions of Section~\ref{sec:algorithm}, fix a constant radius $t>0$. After at most
\[
 C_t=\ceilp{\frac{4}{t}\left(s_0-\frac{2Lt}{\alpha}\right)}
\]
iterations, either an optimum has been reached or every subsequent update satisfies
\begin{equation}\label{eq:fixedrec}
 \Delta_{k+1}\le \Delta_k-
 \min\!\left\{\frac{\alpha}{4L}\Delta_{k},\frac{\Delta_{k}^2}{2Lt^2}\right\}.
\end{equation}
The objective values converge to $f_\star$, and $\dist(x_k,\Xs)\to0$. For fixed initialization and parameters, we obtain
\begin{equation}\label{eq:fixedrate}
 \Delta_k=O\!\left(\frac{Lt^2}{k}\right).
\end{equation}
\end{theorem}
This upper bound holds for every positive constant radius, regardless of whether the ball covers $\X$.

\begin{proposition}[No uniform fixed-radius geometric bound]\label{prop:lower}
For every fixed $t>0$, there is a smooth strongly convex quadratic on a compact convex set of diameter less than $t$ for which no constants $A<\infty$ and $q_{\mathrm{geom}}\in(0,1)$ give
\[
 \Delta_k\le Aq_{\mathrm{geom}}^k
\]
uniformly over all feasible initial points and all $k$.
\end{proposition}
The proof rescales the fixed ball-constrained instance of \citet{halbey2026} and applies Proposition~\ref{prop:fw}. Since the strongly convex counterexample also satisfies quadratic growth, the proposition rules out a fixed-radius geometric bound uniform over all initial points under our main assumptions.

\section{Conclusion}
\begingroup\color{black}
An accepted short step can control the next local search radius without knowing the remaining objective error. For the resulting short-step variant of \algname{BroxGD}, we proved linear convergence under quadratic growth on nonempty, closed, convex sets. The radius invariant explains \termrev{why the adaptive radius rule works}; the contraction-or-halving argument counts every exact local-oracle call, including initialization. A gradient-based initial radius removes that phase, and \textcolor{black}{an $R$-linear convergence bound for the objective values} yields convergence of the iterates to a single solution. These guarantees need a known smoothness bound but no growth constant, optimal value, or minimizer.

The reductions to \algname{ProjGD} and short-step \algname{FW} explain how radius choices connect familiar methods. The fixed-radius analysis shows why adaptation matters: a constant radius has a sublinear upper bound, and in general no geometric bound holds uniformly over all initial points. The present guarantee also differs from an accepted-step rate for a backtracking method, since there are no rejected local solves to account for.

The established computational resource is an exact local-oracle call. We do not infer a runtime advantage over projection or global linear minimization from this count. Efficient implementations for specific feasible sets and guarantees for inexact local solves remain important next questions. The latter must account for the use of exact optimality in both the local comparison \termrev{and the objective-error bound after a full interior step}.
\endgroup

\begingroup\color{black}
\section*{Acknowledgements}
This work was supported by funding from King Abdullah University of Science and Technology (KAUST): i) KAUST Baseline Research Scheme, ii) Center of Excellence for Generative AI (award no. 5940).
\endgroup

\bibliographystyle{plainnat}
\bibliography{references}
\clearpage
\appendix
\etocdepthtag.toc{appendix}
\section{Proofs}\label{app:proofs}

\subsection{Well-posedness and descent}
\begin{proof}[\normalfont\bfseries Proof of Lemma~\ref{lem:descent}]
The local feasible set $\X\cap\B(x_k,t_k)$ contains $x_k$ and is closed and bounded, hence compact. The continuous linear objective attains its minimum on this set. Comparing any minimizer $z_k$ with $x_k$ gives
\[
 \inner{\nabla f(x_k)}{z_k}\overset{\textcolor{black}{\eqref{eq:oracle}}}{\textcolor{black}{\le}}\inner{\nabla f(x_k)}{x_k},
 \qquad \inner{\nabla f(x_k)}{x_k-z_k}\ge0.
\]
We first show that equality characterizes optimality on all of $\X$. Suppose the local gap is zero and fix $y\in\X\setminus\left\{x_k\right\}$. Choose
\[
 \theta=\min\!\left\{1,\frac{t_k}{\|y-x_k\|}\right\},
 \qquad w=x_k+\theta(y-x_k).
\]
Since $0<\theta\le1$, convexity gives $w\in\X$, and $\|w-x_k\|=\theta\|y-x_k\|\le t_k$. Thus $w$ is feasible for the local oracle. Optimality of $z_k$ implies
\[
 0=\inner{\nabla f(x_k)}{x_k-z_k}
 \overset{\textcolor{black}{\eqref{eq:oracle}}}{\textcolor{black}{\ge}}\inner{\nabla f(x_k)}{x_k-w}
 =\theta\inner{\nabla f(x_k)}{x_k-y}.
\]
Dividing by $\theta>0$ and reversing the displacement yields $\inner{\nabla f(x_k)}{y-x_k}\ge0$. The same inequality holds for $y=x_k$. Convexity of $f$ now gives
\[
 f(y)\ge f(x_k)+\inner{\nabla f(x_k)}{y-x_k}\ge f(x_k)
 \qquad(y\in\X),
\]
so $x_k$ is optimal. Conversely, if $x_k$ is optimal, the constrained first-order condition gives $\inner{\nabla f(x_k)}{z_k-x_k}\ge0$. Hence the local gap is nonpositive. Since it is also nonnegative, it must be zero.

At a nonterminal iteration the gap is positive, so $z_k\ne x_k$. The definition of the accepted length gives
\[
 0<\ell_k\overset{\textcolor{black}{\eqref{eq:length}}}{\textcolor{black}{\le}}\|z_k-x_k\|\le t_k,
 \qquad L\ell_k\overset{\textcolor{black}{\eqref{eq:length}}}{\textcolor{black}{\le}}
 \frac{\inner{\nabla f(x_k)}{x_k-z_k}}{\|z_k-x_k\|}.
\]
The update is the convex combination
\[
 x_{k+1}\overset{\textcolor{black}{\eqref{eq:update}}}{\textcolor{black}{=}}\left(1-\frac{\ell_k}{\|z_k-x_k\|}\right)x_k
             +\frac{\ell_k}{\|z_k-x_k\|}z_k,
\]
so $x_{k+1}\in\X$. Also, $\frac{t_k}{2}\le t_{k+1}=\max\left\{\ell_k,\frac{t_k}{2}\right\}\le t_k$. These facts justify all subsequent oracle calls and updates by induction from $x_0\in\X$ and $t_0>0$.

To bound the decrease, substitute $x_{k+1}-x_k=\ell_k\frac{z_k-x_k}{\|z_k-x_k\|}$ into the smoothness inequality. Since this displacement has norm $\ell_k$, we obtain
\[
 \rev{\Delta_k}-\rev{\Delta_{k+1}}\overset{\textcolor{black}{\eqref{eq:smooth},\eqref{eq:update}}}{\textcolor{black}{\ge}}
 \ell_k\frac{\inner{\nabla f(x_k)}{x_k-z_k}}{\|z_k-x_k\|}
 -\frac{L}{2}\ell_k^2.
\]
Multiplying the upper bound on $L\ell_k$ by $\frac{\ell_k}{2}$ gives
\[
 \frac{L}{2}\ell_k^2\overset{\textcolor{black}{\eqref{eq:length}}}{\textcolor{black}{\le}}
 \frac{\ell_k}{2}\frac{\inner{\nabla f(x_k)}{x_k-z_k}}{\|z_k-x_k\|}.
\]
Subtracting this from the linear term proves
\[
 \rev{\Delta_k}-\rev{\Delta_{k+1}}\overset{\textcolor{black}{\eqref{eq:smooth},\eqref{eq:update},\eqref{eq:length}}}{\textcolor{black}{\ge}}
 \frac{\ell_k}{2}\frac{\inner{\nabla f(x_k)}{x_k-z_k}}{\|z_k-x_k\|}>0.
\]
The last inequality uses positivity of both the accepted length and the local gap.

For a short step, the first entry of the minimum defining $\ell_k$ is active, so
\[
 \ell_k\overset{\textcolor{black}{\eqref{eq:length}}}{\textcolor{black}{=}}-\frac{\inner{\nabla f(x_k)}{z_k-x_k}}{L\|z_k-x_k\|}
 =\frac{\inner{\nabla f(x_k)}{x_k-z_k}}{L\|z_k-x_k\|}.
\]
Substituting this expression into both terms of the model decrease gives
\[
\begin{relchain}
\rev{\Delta_k}-\rev{\Delta_{k+1}} & \overset{\textcolor{black}{\eqref{eq:smooth},\eqref{eq:update},\eqref{eq:length}}}{\textcolor{black}{\ge}} & \frac{\inner{\nabla f(x_k)}{x_k-z_k}^2}{L\|z_k-x_k\|^2}
      -\frac{\inner{\nabla f(x_k)}{x_k-z_k}^2}{2L\|z_k-x_k\|^2}\\[4pt]
 & = & \frac{\inner{\nabla f(x_k)}{x_k-z_k}^2}{2L\|z_k-x_k\|^2}
 \overset{\textcolor{black}{\eqref{eq:length}}}{\textcolor{black}{=}}\frac{L}{2}\ell_k^2.
\end{relchain}
\]
This is the short-step estimate.

For a full step, $\ell_k=\|z_k-x_k\|$. The minimum defining $\ell_k$ can select this length only if
\[
 \|z_k-x_k\|\overset{\textcolor{black}{\eqref{eq:length}}}{\textcolor{black}{\le}}
 \frac{\inner{\nabla f(x_k)}{x_k-z_k}}{L\|z_k-x_k\|},
\]
which implies $L\|z_k-x_k\|^2\le\inner{\nabla f(x_k)}{x_k-z_k}$. Substituting the full length into the model and then using this inequality gives
\[
\begin{relchain}
\rev{\Delta_k}-\rev{\Delta_{k+1}} & \overset{\textcolor{black}{\eqref{eq:smooth},\eqref{eq:update}}}{\textcolor{black}{\ge}} & \inner{\nabla f(x_k)}{x_k-z_k}-\frac{L}{2}\|z_k-x_k\|^2\\[4pt]
 & \overset{\textcolor{black}{\eqref{eq:length}}}{\textcolor{black}{\ge}} & \inner{\nabla f(x_k)}{x_k-z_k}
       -\frac{1}{2}\inner{\nabla f(x_k)}{x_k-z_k}\\[4pt]
 & = & \frac{1}{2}\inner{\nabla f(x_k)}{x_k-z_k}.
\end{relchain}
\]
Equality in the length test is covered by this full-step argument.
\end{proof}

\subsection{Scaling and interior steps}
\begin{proof}[\normalfont\bfseries Proof of Lemma~\ref{lem:scaling}]
Fix $y\in\X\setminus\left\{x_k\right\}$. We compare the oracle solution with a point on the segment from $x_k$ to $y$ that lies in the ball. Set
\[
 \theta=\min\!\left\{1,\frac{t_k}{\|y-x_k\|}\right\},
 \qquad w=x_k+\theta(y-x_k).
\]
Convexity gives $w\in\X$, and $\|w-x_k\|=\theta\|y-x_k\|\le t_k$. Since $w$ is feasible for the local oracle,
\[
 \inner{\nabla f(x_k)}{x_k-z_k}
 \overset{\textcolor{black}{\eqref{eq:oracle}}}{\textcolor{black}{\ge}}\inner{\nabla f(x_k)}{x_k-w}
 =\theta\inner{\nabla f(x_k)}{x_k-y}.
\]
This proves the scaling inequality without using convexity of $f$.

Now suppose $x_k$ is nonoptimal and choose a nearest solution $x_{\star,k}=\Proj_{\Xs}(x_k)$. It exists because $\Xs$ is nonempty, closed, and convex, and its distance from $x_k$ is positive. Convexity of $f$ gives
\[
 f_\star\ge f(x_k)+\inner{\nabla f(x_k)}{x_{\star,k}-x_k},
\]
so the scaling inequality with $y=x_{\star,k}$ yields
\[
 \inner{\nabla f(x_k)}{x_k-z_k}
 \overset{\textcolor{black}{\eqref{eq:scaling}}}{\textcolor{black}{\ge}}\min\!\left\{1,\frac{t_k}{\dist(x_k,\Xs)}\right\}\rev{\Delta_k}.
\]
The right-hand side is positive. We may therefore divide by $\|z_k-x_k\|>0$ and use $\|z_k-x_k\|\le t_k$ to obtain
\[
\begin{relchain}
\frac{\inner{\nabla f(x_k)}{x_k-z_k}}{\|z_k-x_k\|} & \overset{\textcolor{black}{\eqref{eq:scaling},\eqref{eq:feasible}}}{\textcolor{black}{\ge}} & \frac{\rev{\Delta_k}}{t_k}
       \min\!\left\{1,\frac{t_k}{\dist(x_k,\Xs)}\right\}\\[4pt]
 & = & \frac{\rev{\Delta_k}}{\max\left\{t_k,\dist(x_k,\Xs)\right\}}.
\end{relchain}
\]
For the last equality, if $t_k\le\dist(x_k,\Xs)$ the minimum is $\frac{t_k}{\dist(x_k,\Xs)}$, while if $t_k\ge\dist(x_k,\Xs)$ it is one. These two cases give the displayed maximum in the denominator.
\end{proof}

\begin{proof}[\normalfont\bfseries Proof of Lemma~\ref{lem:reset}]
Since $\|z_k-x_k\|<t_k$, there is room to move a positive distance from $z_k$ toward any other feasible point. Fix $y\in\X\setminus\left\{z_k\right\}$ and take
\[
 \theta=\min\!\left\{1,\frac{t_k-\|z_k-x_k\|}{\|y-z_k\|}\right\}>0,
 \qquad w=(1-\theta)z_k+\theta y.
\]
Then $w\in\X$ by convexity and
\[
 \|w-x_k\|\le\|z_k-x_k\|+\theta\|y-z_k\|\le t_k.
\]
Oracle optimality gives $\inner{\nabla f(x_k)}{w-z_k}\ge0$. Substituting $w-z_k=\theta(y-z_k)$ and dividing by $\theta>0$ shows that
\[
 \inner{\nabla f(x_k)}{y-z_k}\ge0\qquad(y\in\X).
\]
The case $y=z_k$ holds with equality. Thus $z_k$ minimizes the linear objective over all of $\X$, irrespective of convexity of $f$.

For any solution $x_\star$, global linear optimality and convexity of $f$ give, respectively,
\[
 \inner{\nabla f(x_k)}{z_k-x_k}
 \le\inner{\nabla f(x_k)}{x_\star-x_k}
 \le f_\star-f(x_k).
\]
Because the step is full, $x_{k+1}=z_k$. Applying smoothness to this point and then the preceding bound gives
\[
\begin{relchain}
f(x_{k+1}) & \overset{\textcolor{black}{\eqref{eq:smooth},\eqref{eq:update}}}{\textcolor{black}{\le}} & f(x_k)+\inner{\nabla f(x_k)}{z_k-x_k}
                 +\frac{L}{2}\|z_k-x_k\|^2\\[4pt]
 & \le & f_\star+\frac{L}{2}\|z_k-x_k\|^2.
\end{relchain}
\]
Subtract $f_\star$, multiply by $\frac{2}{L}$, and take nonnegative square roots. Since $t_{k+1}\ge\ell_k=\|z_k-x_k\|$, we obtain
\[
 t_{k+1}\overset{\textcolor{black}{\eqref{eq:update}}}{\textcolor{black}{\ge}}\|z_k-x_k\|\ge
 \sqrt{\frac{2\rev{\Delta_{k+1}}}{L}}.
\]
\end{proof}

\begin{revision}
\subsection{The radius invariant and initialization}\label{app:invariant}
We use $\rev{\Delta_k},s_k$ from \eqref{eq:analysisnotation} and $G_k,D_k$ from \eqref{eq:stepnotation}. This makes the radius scale visible throughout the case analysis.
\begin{proof}[\normalfont\bfseries Proof of Lemma~\ref{lem:envelope}]
First consider a nonterminal iterate with $t_k\ge\alpha s_k/(2L)$. Descent implies $s_{k+1}\le s_k$. A full boundary step preserves the invariant because
\[
 t_{k+1}\overset{\eqref{eq:update}}{=}t_k
 \overset{\eqref{eq:invariant}}{\ge}\frac{\alpha s_k}{2L}
 \overset{\eqref{eq:descent},\eqref{eq:analysisnotation}}{\ge}\frac{\alpha s_{k+1}}{2L}.
\]
For a full interior step, \termrev{the objective-error estimate for a full interior step} gives
\[
 t_{k+1}\overset{\eqref{eq:reset},\eqref{eq:analysisnotation}}{\ge}
 \sqrt{\frac{\alpha}{L}}\,s_{k+1}
 \overset{\eqref{eq:scales}}{\ge}\frac{\alpha s_{k+1}}{2L}.
\]
The last comparison uses $0<\alpha/L\le1$.

\auditrev{For the remaining case, first note that whenever $t_k\le s_k$, both $t_k$ and $\dist(x_k,\Xs)$ are at most $s_k$. Irrespective of the step type, this gives}
\begin{equation}\label{eq:small-slope}
 \frac{G_k}{D_k}
 \overset{\eqref{eq:slopes},\eqref{eq:stepnotation}}{\ge}
 \frac{\rev{\Delta_k}}{\max\left\{t_k,\dist(x_k,\Xs)\right\}}
 \overset{\eqref{eq:qg},\eqref{eq:scales},\eqref{eq:analysisnotation}}{\ge}
 \frac{\rev{\Delta_k}}{s_k}
 \overset{\eqref{eq:analysisnotation}}{=}\frac{\alpha s_k}{2}.
\end{equation}
Here $s_k>0$ because the iterate is nonoptimal. \auditrev{If the step is short and $t_k\le s_k$, its length therefore yields}
\[
 t_{k+1}\overset{\eqref{eq:update}}{\ge}\ell_k
 \overset{\eqref{eq:length},\eqref{eq:small-slope}}{\ge}\frac{\alpha s_k}{2L}
 \overset{\eqref{eq:descent},\eqref{eq:analysisnotation}}{\ge}\frac{\alpha s_{k+1}}{2L}.
\]
If instead $t_k>s_k$, the safeguard gives
\[
 t_{k+1}\overset{\eqref{eq:update}}{\ge}\frac{t_k}{2}
 >\frac{s_k}{2}
 \overset{\eqref{eq:scales}}{\ge}\frac{\alpha s_k}{2L}
 \overset{\eqref{eq:descent},\eqref{eq:analysisnotation}}{\ge}\frac{\alpha s_{k+1}}{2L}.
\]
This covers every step type, including a step reaching an optimum, when $s_{k+1}=0$.

Now suppose the invariant fails. Since $t_k<\alpha s_k/(2L)\le s_k$, \eqref{eq:small-slope} applies and gives
\[
 \frac{G_k}{LD_k}\overset{\eqref{eq:small-slope}}{\ge}\frac{\alpha s_k}{2L}
 >t_k\overset{\eqref{eq:feasible},\eqref{eq:stepnotation}}{\ge}D_k.
\]
Thus clipping selects a full step. An interior step establishes the invariant by \termrev{the argument for a full interior step above}. Otherwise $D_k=\ell_k=t_k$ and the radius remains unchanged. Writing $\rev{\delta_k}=\rev{\Delta_k}-\rev{\Delta_{k+1}}$, the model estimate gives
\begin{equation}\label{eq:initial-decrease}
 \rev{\delta_k}\overset{\eqref{eq:descent},\eqref{eq:small-slope},\eqref{eq:analysisnotation}}{\ge}
 \frac{\alpha t_ks_k}{2}-\frac{Lt_k^2}{2}
 >\frac{\alpha t_ks_k}{4},
\end{equation}
where the strict inequality uses $Lt_k<\alpha s_k/2$. The difference-of-squares identity then yields
\[
\begin{relchain}
 s_k-s_{k+1}
 &\overset{\eqref{eq:analysisnotation}}{=}&
 \frac{2\rev{\delta_k}}{\alpha(s_k+s_{k+1})}\\[3pt]
 &\overset{\eqref{eq:descent},\eqref{eq:analysisnotation}}{\ge}&
 \frac{\rev{\delta_k}}{\alpha s_k}
 \overset{\eqref{eq:initial-decrease}}{>}\frac{t_k}{4}.
\end{relchain}
\]
Until the invariant holds, boundary steps keep the radius at $t_0$ and reduce $s_k$ by at least $t_0/4$ per update. An interior step ends the phase. Summing until $s_k\le2Lt_0/\alpha$ gives the ceiling bound $N\le C_0$ in \eqref{eq:catchup}; Appendix~\ref{app:count} explicitly accounts for a possible final interior step.
\end{proof}

\subsection{Contraction or halving}
The invariant supplies a useful radius; the remaining question is whether the accepted step reduces the error or the radius.
\begin{proof}[\normalfont\bfseries Proof of Lemma~\ref{lem:dichotomy}]
Write $\rev{\delta_k}=\rev{\Delta_k}-\rev{\Delta_{k+1}}$. We prove $\rev{\delta_k}\ge\alpha \rev{\Delta_k}/(4L)$ except in the stated halving case.

Suppose first that $\alpha s_k/(2L)\le t_k\le s_k$. Quadratic growth and the invariant give
\[
 \min\left\{1,\frac{t_k}{\dist(x_k,\Xs)}\right\}
 \overset{\eqref{eq:qg},\eqref{eq:scales},\eqref{eq:analysisnotation}}{\ge}
 \frac{t_k}{s_k}
 \overset{\eqref{eq:invariant}}{\ge}\frac{\alpha}{2L}.
\]
The local-gap comparison yields $G_k\ge\alpha \rev{\Delta_k}/(2L)$. For a full step,
\[
 \rev{\delta_k}\overset{\eqref{eq:full},\eqref{eq:stepnotation},\eqref{eq:analysisnotation}}{\ge}\frac{G_k}{2}
 \overset{\eqref{eq:slopes},\eqref{eq:invariant},\eqref{eq:qg},\eqref{eq:analysisnotation}}{\ge}\frac{\alpha \rev{\Delta_k}}{4L}.
\]
For a short step, use the \termrev{normalized local linear gap $G_k/D_k$} instead:
\[
 \rev{\delta_k}\overset{\eqref{eq:short},\eqref{eq:stepnotation},\eqref{eq:analysisnotation}}{\ge}
 \frac{1}{2L}\left(\frac{G_k}{D_k}\right)^2
 \overset{\eqref{eq:small-slope}}{\ge}\frac{\alpha^2s_k^2}{8L}
 \overset{\eqref{eq:analysisnotation}}{=}\frac{\alpha \rev{\Delta_k}}{4L}.
\]
Hence every step in this radius range contracts the error.

Now suppose $t_k>s_k$. A nearest solution lies in the local ball, so $G_k\ge \rev{\Delta_k}$ by \eqref{eq:slopes}. A full step therefore gives
\[
 \rev{\delta_k}\overset{\eqref{eq:full},\eqref{eq:slopes},\eqref{eq:analysisnotation}}{\ge}\frac{\rev{\Delta_k}}{2}
 \overset{\eqref{eq:scales}}{\ge}\frac{\alpha \rev{\Delta_k}}{4L}.
\]
For a short step with $\ell_k>t_k/2$, we have
\[
\begin{relchain}
 \rev{\delta_k}&\overset{\eqref{eq:length},\eqref{eq:short},\eqref{eq:analysisnotation}}{\ge}&
 \frac L2\ell_k^2>\frac L8t_k^2>\frac L8s_k^2\\[3pt]
 &\overset{\eqref{eq:analysisnotation}}{=}&\frac{L\rev{\Delta_k}}{4\alpha}
 \overset{\eqref{eq:scales}}{\ge}\frac{\alpha \rev{\Delta_k}}{4L}.
\end{relchain}
\]
\termrev{Thus failure to satisfy the contraction bound forces} a short step with $t_k>s_k$ and $\ell_k\le t_k/2$. In that case
\[
 t_{k+1}\overset{\eqref{eq:update}}{=}\max\left\{\ell_k,\frac{t_k}{2}\right\}=\frac{t_k}{2}.
\]
Equality $\ell_k=t_k/2$ is included, so the cases exhaust all updates.
\end{proof}
Even in the halving case, the objective strictly decreases. Using $G_k\ge \rev{\Delta_k}$ and $D_k\le t_k$ gives
\[
 \rev{\delta_k}:=\rev{\Delta_k}-\rev{\Delta_{k+1}}
 \overset{\eqref{eq:short},\eqref{eq:stepnotation},\eqref{eq:analysisnotation}}{\ge}\frac{G_k^2}{2LD_k^2}
 \overset{\eqref{eq:slopes},\eqref{eq:feasible},\eqref{eq:analysisnotation}}{\ge}\frac{\rev{\Delta_k}^2}{2Lt_k^2}>0.
\]
The relative decrease $\rev{\Delta_k}/(2Lt_k^2)$ need not be bounded away from zero when the radius is excessive. This is why the proof counts radius reductions in addition to objective contractions.
\end{revision}

\subsection{Total iteration count}
\label{app:count}
We count the initialization steps, objective contractions, and radius halvings separately, retaining the ceiling terms to obtain a sharper bound. Let
\begin{equation}\label{eq:q}
 q=1-\frac{\alpha}{4L},
\end{equation}
and let $C_0$ be as in \eqref{eq:catchup}. The sharper form of Theorem~\ref{thm:main} is
\begin{equation}\label{eq:complexity}
 K(\varepsilon)\le
 \underbrace{\ceilp{\frac4{t_0}\left(\rev{s_0}-\frac{2Lt_0}{\alpha}\right)}}_{\text{initialization}}
 +\underbrace{\left\lceil\frac{\log\!\left(\frac{\rev{\Delta_0}}{\varepsilon}\right)}{-\log\!\left(1-\frac{\alpha}{4L}\right)}\right\rceil}_{\text{objective contractions}}
 +\underbrace{\ceilp{\log_2\frac{t_0\sqrt{\alpha}}{\sqrt{2\varepsilon}}}}_{\text{radius halvings}}.
\end{equation}

To count these updates, denote the second and third terms of \eqref{eq:complexity} by
\begin{align}
 P_\varepsilon&=\left\lceil\frac{\log\!\left(\frac{\rev{\Delta_0}}{\varepsilon}\right)}{-\log\!\left(1-\frac{\alpha}{4L}\right)}\right\rceil,\label{eq:productivecount}\\
 H_\varepsilon&=\ceilp{\log_2\frac{t_0\sqrt{\alpha}}{\sqrt{2\varepsilon}}}.\label{eq:halvingcount}
\end{align}
\begin{proof}[\normalfont\bfseries Proof of Theorem~\ref{thm:main}]
If $x_0$ is optimal, the zero-gap characterization proves the stopping assertion. Otherwise fix $0<\varepsilon<\rev{\Delta_0}$. We first bound the initialization phase, then count the updates after the invariant holds.

If the invariant holds at $x_0$, take $N=0$. Otherwise every step before it holds is a full boundary step, except possibly for a final full interior step. After $j$ boundary steps, their decreases telescope to give
\[
 \rev{s_j}
 \overset{\textcolor{black}{\eqref{eq:catchupstep}}}{\textcolor{black}{\le}}\rev{s_0}-\frac{jt_0}{4}.
\]
While these steps keep the radius equal to $t_0$, the invariant holds as soon as the left-hand side is at most $\frac{2Lt_0}{\alpha}$. A sufficient condition is therefore
\begin{align*}
 \rev{s_0}-\frac{jt_0}{4}
 &\le\frac{2Lt_0}{\alpha},\\
 j&\ge\frac{4}{t_0}\left(
 \rev{s_0}-\frac{2Lt_0}{\alpha}\right).
\end{align*}
Taking the ceiling of the positive part gives $C_0$. If no interior step occurs first, $C_0$ boundary steps suffice. If an interior step occurs at an earlier index $j<C_0$, it establishes the invariant at $j+1\le C_0$. Thus $N\le C_0$. If the target accuracy or an optimum is reached during this phase, the claimed bound already holds.

After index $N$, call an iteration \emph{contracting} if $\rev{\Delta_{k+1}}\le q\rev{\Delta_k}$. \termrev{All other iterations fail to satisfy this contraction bound.} They still strictly decrease the objective, but do not attain this particular contraction factor. The two classes are disjoint, and every \termrev{iteration not satisfying \eqref{eq:contract}} halves the radius. A contracting iteration may also halve it, which can only improve the radius bound used below.

After $m$ contracting iterations, objective decrease in the intervening iterations gives
\[
 \rev{\Delta_k}\overset{\textcolor{black}{\eqref{eq:contract},\eqref{eq:descent},\eqref{eq:q}}}{\textcolor{black}{\le}} q^m\rev{\Delta_N}
 \overset{\textcolor{black}{\eqref{eq:descent}}}{\textcolor{black}{\le}} q^m\rev{\Delta_0}.
\]
Since $0<q<1$, the last expression is at most $\varepsilon$ whenever
\[
 m\ge\frac{\log\!\left(\frac{\rev{\Delta_0}}{\varepsilon}\right)}{-\log q}.
\]
Hence $P_\varepsilon$ contracting iterations suffice. In particular, there can be at most $P_\varepsilon$ such updates up to and including the first update that reaches the target.

To count \termrev{iterations not satisfying \eqref{eq:contract}}, note that the radius invariant holds after the initialization phase, so Lemma~\ref{lem:dichotomy} applies. Suppose the $H$-th \termrev{iteration not satisfying \eqref{eq:contract}} starts at index $k$ before the target has been reached. Since this iteration does not satisfy the contraction bound and $\rev{\Delta_k}>\varepsilon$, the lemma gives
\[
 t_k>\rev{s_k}
 >\sqrt{\frac{2\varepsilon}{\alpha}}.
\]
The preceding $H-1$ \termrev{iterations not satisfying \eqref{eq:contract}} have each halved the radius. Since no iteration increases it,
\[
 t_k\overset{\textcolor{black}{\eqref{eq:feasible},\eqref{eq:halve}}}{\textcolor{black}{\le}}\frac{t_0}{2^{H-1}}.
\]
Combining the two bounds and rearranging yields
\[
 2^{H-1}<\frac{t_0}{\sqrt{\frac{2\varepsilon}{\alpha}}}
 =\frac{t_0\sqrt{\alpha}}{\sqrt{2\varepsilon}},
 \qquad
 H-1<\log_2\frac{t_0\sqrt{\alpha}}{\sqrt{2\varepsilon}}.
\]
If the logarithm's argument is at most one, the logarithm is nonpositive, which contradicts $H-1\ge0$. Thus no \termrev{iteration not satisfying \eqref{eq:contract}} can start before the target in this case. If the argument exceeds one, the strict inequality and integrality of $H$ imply
\[
 H\le\left\lceil\log_2\frac{t_0\sqrt{\alpha}}{\sqrt{2\varepsilon}}\right\rceil.
\]
If $H$ were larger than this ceiling, then $H-1$ would be at least the ceiling and hence at least the logarithm, contradicting the strict inequality. Both cases are therefore covered by $H\le H_\varepsilon$.

Every update before the target belongs to the initialization phase or one of the two classes. The preceding finite bounds also show that the target must be reached. Their sum gives
\[
 K(\varepsilon)\overset{\textcolor{black}{\eqref{eq:catchup},\eqref{eq:productivecount},\eqref{eq:halvingcount}}}{\textcolor{black}{\le}} C_0+P_\varepsilon+H_\varepsilon,
\]
which is the sharper count stated above.

For the bound in the main text, replace each ceiling by its argument plus one. Also, $-\log(1-u)\ge u$ for $0<u<1$, so
\[
 \frac{1}{-\log\!\left(1-\frac{\alpha}{4L}\right)}
 \le\frac{4L}{\alpha}.
\]
Substituting these bounds gives the three terms in \eqref{eq:nonasymptotic} and the additive constant three. With initialization and problem constants fixed, $C_0$ is constant, $P_\varepsilon=O\!\left(\frac{L}{\alpha}\log\frac{1}{\varepsilon}\right)$, and $H_\varepsilon=O\!\left(\log\frac{1}{\varepsilon}\right)$. Since $\frac{L}{\alpha}\ge1$, their sum has the asserted asymptotic order. Each counted update uses one oracle call. If an optimum is reached, the zero-gap test there may require one additional call to certify termination.
\end{proof}

\subsection{A geometric bound and point convergence}\label{app:fullsequence}
The same argument yields a geometric bound on the entire objective-error sequence and convergence of the iterates to a single solution.
\begin{corollary}[$R$-linear objective and point convergence]\label{cor:rlinear}
\auditrev{Under the assumptions of Theorem~\ref{thm:main},} let $N$ satisfy
\[
 t_N\ge\frac{\alpha}{2L}\rev{s_N}.
\]
If $\rev{\Delta_N}>0$, define
\[
 \simplerev{B_N=\frac{2L^2t_N^2}{\alpha}.}
\]
Then, for every $n\ge0$,
\begin{equation}\label{eq:rlinear}
 \simplerev{\Delta_{N+n}\le B_N\exp\!\left(-\frac{\alpha n}{6L}\right).}
\end{equation}
The trajectory has finite length and converges to one $x_\infty\in\Xs$, with
\[
 \simplerev{\|x_{N+n}-x_\infty\|=O\!\left(\exp\!\left(-\frac{\alpha n}{12L}\right)\right).}
\]
\simplerev{Such an index $N$ can be chosen to satisfy \eqref{eq:catchup}}, and $B_N\le\frac{2L^2t_0^2}{\alpha}$. If an optimum is reached, keep all subsequent iterates fixed.
\end{corollary}
\begin{proof}[\normalfont\bfseries Proof of Corollary~\ref{cor:rlinear}]
\simplerev{We retain the sharper constants in the proof. Set $a=-\log q$, $b=\log4$, and $\lambda=ab/(a+b)$. We first establish}
\begin{equation}\label{eq:rlinear-exact}
 \simplerev{\Delta_{N+n}\le\Delta_N^{\frac{b}{a+b}}B_N^{\frac{a}{a+b}}e^{-\lambda n}.}
\end{equation}

First consider $n$ updates starting at $N$ before an optimum is reached. Let $P$ and $H$ count \termrev{iterations satisfying the contraction bound and those that do not} as in the preceding proof. Then $P+H=n$. Each contracting iteration reduces the error by a factor at most $q$, and the other iterations do not increase it. Thus
\begingroup\makeatletter\def\tagform@#1{\maketag@@@{\textcolor{black}{(#1)}}}\makeatother
\begin{equation}\label{eq:objective-count-bound}
\rev{\Delta_{N+n}}\overset{\textcolor{black}{\eqref{eq:contract},\eqref{eq:descent},\eqref{eq:q}}}{\textcolor{black}{\le}}\rev{\Delta_N}q^P
 =\rev{\Delta_N}e^{-aP}.
\end{equation}
\endgroup
Each \termrev{iteration not satisfying \eqref{eq:contract}} halves the radius, so $t_{N+n}\le t_N2^{-H}$. Squaring the invariant and solving for the objective error gives
\[
 t_{N+n}^2\overset{\textcolor{black}{\eqref{eq:invariant}}}{\textcolor{black}{\ge}}\frac{\alpha}{2L^2}\rev{\Delta_{N+n}},
\]
and hence
\begingroup\makeatletter\def\tagform@#1{\maketag@@@{\textcolor{black}{(#1)}}}\makeatother
\begin{equation}\label{eq:radius-count-bound}
\rev{\Delta_{N+n}}\overset{\textcolor{black}{\eqref{eq:invariant}}}{\textcolor{black}{\le}}\frac{2L^2}{\alpha}t_{N+n}^2
 \overset{\textcolor{black}{\eqref{eq:feasible},\eqref{eq:halve}}}{\textcolor{black}{\le}}\frac{2L^2t_N^2}{\alpha}4^{-H}=B_Ne^{-bH}.
\end{equation}
\endgroup
To combine these bounds, raise the first to $\frac{b}{a+b}$ and the second to $\frac{a}{a+b}$. The exponents are positive and sum to one. Multiplying the resulting inequalities therefore gives
\[
\begin{relchain}
\rev{\Delta_{N+n}} & = & \rev{\Delta_{N+n}}^{\frac{b}{a+b}}
   \rev{\Delta_{N+n}}^{\frac{a}{a+b}}\\[4pt]
 & \overset{\textcolor{black}{\eqref{eq:objective-count-bound},\eqref{eq:radius-count-bound}}}{\textcolor{black}{\le}} & \rev{\Delta_N}^{\frac{b}{a+b}}B_N^{\frac{a}{a+b}}
    \exp\!\left(-\frac{abP+abH}{a+b}\right)\\[4pt]
 & = & \rev{\Delta_N}^{\frac{b}{a+b}}B_N^{\frac{a}{a+b}}e^{-\lambda n}.
\end{relchain}
\]
This proves the geometric bound. If an optimum is reached, all later objective errors are zero, so the same bound continues to hold when the iterate is kept fixed.

We next turn objective convergence into convergence of the iterates. Put $\auditrev{A_N}=\rev{\Delta_N}^{\frac{b}{a+b}}B_N^{\frac{a}{a+b}}$. The accepted-length formula implies that the local slope is at least $L\ell_k$. Applying this to the descent estimate gives
\[
 \rev{\Delta_k}-\rev{\Delta_{k+1}}
 \overset{\textcolor{black}{\eqref{eq:descent}}}{\textcolor{black}{\ge}}\frac{\ell_k}{2}\frac{\inner{\nabla f(x_k)}{x_k-z_k}}{\|z_k-x_k\|}
 \overset{\textcolor{black}{\eqref{eq:length}}}{\textcolor{black}{\ge}}\frac{L}{2}\ell_k^2.
\]
Since $f(x_{k+1})\ge f_\star$, it follows that
\begingroup\makeatletter\def\tagform@#1{\maketag@@@{\textcolor{black}{(#1)}}}\makeatother
\begin{equation}\label{eq:summable-length-bound}
\ell_{N+j}^2\overset{\textcolor{black}{\eqref{eq:descent},\eqref{eq:length}}}{\textcolor{black}{\le}}\frac{2}{L}\rev{\Delta_{N+j}}
 \overset{\textcolor{black}{\eqref{eq:rlinear-exact}}}{\textcolor{black}{\le}}\frac{2\auditrev{A_N}}{L}e^{-\lambda j},
 \qquad
 \ell_{N+j}\overset{\textcolor{black}{\eqref{eq:descent},\eqref{eq:length},\eqref{eq:rlinear-exact}}}{\textcolor{black}{\le}}\sqrt{\frac{2\auditrev{A_N}}{L}}e^{-\frac{\lambda j}{2}}.
\end{equation}
\endgroup
After termination set $\ell_k=0$, so this bound remains valid. Because $\lambda>0$, the step-length series is summable. For $m>n$, the triangle inequality gives the explicit tail estimate
\[
\begin{relchain}
\|x_{N+m}-x_{N+n}\| & \overset{\textcolor{black}{\eqref{eq:update}}}{\textcolor{black}{\le}} & \sum_{j=n}^{m-1}\ell_{N+j}\\[4pt]
 & \overset{\textcolor{black}{\eqref{eq:summable-length-bound}}}{\textcolor{black}{\le}} & \sqrt{\frac{2\auditrev{A_N}}{L}}\sum_{j=n}^{\infty}e^{-\frac{\lambda j}{2}}
 =\frac{\sqrt{\frac{2\auditrev{A_N}}{L}}}{1-e^{-\frac{\lambda}{2}}}e^{-\frac{\lambda n}{2}}.
\end{relchain}
\]
The right-hand side tends to zero as $n\to\infty$, so the iterates form a Cauchy sequence. Their limit $x_\infty$ belongs to $\X$ because $\X$ is closed. Continuity of $f$ and the geometric error bound imply $f(x_\infty)=f_\star$. Letting $m\to\infty$ in the last inequality proves
\begin{equation}\label{eq:pointtail}
 \|x_{N+n}-x_\infty\|\overset{\textcolor{black}{\eqref{eq:update},\eqref{eq:summable-length-bound}}}{\textcolor{black}{\le}}
 \frac{\sqrt{\frac{2\auditrev{A_N}}{L}}}{1-e^{-\frac{\lambda}{2}}}e^{-\frac{\lambda n}{2}}.
\end{equation}
The finite prefix before $N$ adds only finitely many step lengths, so the entire trajectory has finite length.

Finally, we verify the order of $\lambda$. For $0<\frac{\alpha}{4L}\le\frac{1}{4}$, the elementary logarithm bounds give
\[
 \frac{\alpha}{4L}\le a\overset{\textcolor{black}{\eqref{eq:q}}}{\textcolor{black}{=}}-\log\!\left(1-\frac{\alpha}{4L}\right)
 \le\frac{\frac{\alpha}{4L}}{1-\frac{\alpha}{4L}}
 \overset{\textcolor{black}{\eqref{eq:scales}}}{\textcolor{black}{\le}}\frac{\alpha}{3L}\overset{\textcolor{black}{\eqref{eq:scales}}}{\textcolor{black}{\le}}\frac{1}{3}.
\]
Since $b=\log4$ is fixed,
\[
 \frac{b}{b+\frac{1}{3}}\frac{\alpha}{4L}
 \le\frac{ab}{a+b}=\lambda\le a\le\frac{\alpha}{3L}.
\]
Thus $\lambda=\Theta\left(\frac{\alpha}{L}\right)$. \simplerev{Since $b=\log4\ge1$, the preceding lower bound gives $\lambda\ge3\alpha/(16L)\ge\alpha/(6L)$. Moreover, $\Delta_N\le B_N$ by the invariant, so $A_N\le B_N$. Applying these comparisons to \eqref{eq:rlinear-exact} and \eqref{eq:pointtail} proves the simpler objective and point bounds in the statement.}  \textcolor{black}{Initialization guarantees the existence of an eligible index $N$ satisfying \eqref{eq:catchup}; it does not bound every eligible index.} \textcolor{black}{Moreover,} $B_N\le\frac{2L^2t_0^2}{\alpha}$ follows from $t_N\le t_0$.
\end{proof}

\begingroup
\color{black}
\hypersetup{linkcolor=black}
\subsection{A shorter Lyapunov proof}\label{app:lyapunov}
The radius invariant also yields a direct geometric bound through a single potential. This argument uses the same assumptions as Corollary~\ref{cor:rlinear} and retains Lemma~\ref{lem:envelope}; it replaces the subsequent counting of \termrev{iterations satisfying the contraction bound and those that do not}.

Choose an index $N$ at which \eqref{eq:invariant} holds, with $\rev{\Delta_N}>0$, and define
\begin{equation}\label{eq:lyap-def}
 \rev{\Delta_k}=f(x_k)-f_\star,\qquad r_k=\frac L2t_k^2,\qquad
 V_k=\rev{\Delta_k}+r_k,\qquad \rev{\delta_k}=\rev{\Delta_k}-\rev{\Delta_{k+1}}.
\end{equation}
\simplerev{In particular, this proof gives the interpretable sufficient number of additional updates}
\[
 \simplerev{\widehat n_\varepsilon=\left\lceil\frac{20L}{3\alpha}\left[\log\!\left(\frac{V_N}{\varepsilon}\right)\right]_+\right\rceil,}
\]
\simplerev{where $V_N=\Delta_N+Lt_N^2/2$. The exact contraction and count are retained below.}

At every nonterminal update, the accepted-length formula and descent estimate give
\begin{equation}\label{eq:lyap-descent}
 \rev{\delta_k}\overset{\eqref{eq:length},\eqref{eq:descent}}{\ge}\frac L2\ell_k^2\ge0.
\end{equation}
Consequently, the radius update satisfies
\begin{equation}\label{eq:lyap-radius}
 r_{k+1}\overset{\eqref{eq:lyap-def},\eqref{eq:update}}{=}
 \max\left\{\frac L2\ell_k^2,\frac{r_k}{4}\right\}
 \overset{\eqref{eq:lyap-descent}}{\le}
 \max\left\{\rev{\delta_k},\frac{r_k}{4}\right\}.
\end{equation}
It follows that
\begin{equation}\label{eq:lyap-decrease}
\begin{array}{rcl}
 V_{k+1}
 &\overset{\eqref{eq:lyap-def},\eqref{eq:lyap-radius}}{\le}&
 \rev{\Delta_k}-\rev{\delta_k}+\max\left\{\rev{\delta_k},r_k/4\right\}\\[3pt]
 &=&\max\left\{\rev{\Delta_k},\rev{\Delta_k}-\rev{\delta_k}+r_k/4\right\}\\[3pt]
 &\overset{\eqref{eq:lyap-descent}}{\le}&\rev{\Delta_k}+r_k/4
 \overset{\eqref{eq:lyap-def}}{=}V_k-3r_k/4.
\end{array}
\end{equation}
For $k\ge N$, the invariant implies
\begin{equation}\label{eq:lyap-control}
 r_k\auditrev{\overset{\eqref{eq:invariant},\eqref{eq:analysisnotation},\eqref{eq:lyap-def}}{\ge}}
 \frac{\alpha}{4L}\rev{\Delta_k}.
\end{equation}
Thus
\begin{equation}\label{eq:lyap-share}
 V_k\overset{\eqref{eq:lyap-def},\eqref{eq:lyap-control}}{\le}
 \left(1+\frac{4L}{\alpha}\right)r_k.
\end{equation}
Writing $q_{\mathrm{Lyap}}=1-\frac{3\alpha}{4(4L+\alpha)}\in(0,1)$, we obtain
\begin{equation}\label{eq:lyap-contract}
 V_{k+1}\overset{\eqref{eq:lyap-decrease},\eqref{eq:lyap-share}}{\le}
 q_{\mathrm{Lyap}} V_k,\qquad k\ge N.
\end{equation}
Induction proves
\begin{equation}\label{eq:lyap-rate}
 \rev{\Delta_{N+n}}\overset{\eqref{eq:lyap-def}}{\le}V_{N+n}
 \overset{\eqref{eq:lyap-contract}}{\le}V_Nq_{\mathrm{Lyap}}^n,\qquad n\ge0.
\end{equation}
If an optimum is reached, extend the sequence for this proof by keeping $x_k$ fixed, setting $\ell_k=0$, and halving $t_k$. Then \eqref{eq:lyap-descent}--\eqref{eq:lyap-contract} remain valid; no further oracle calls are needed for this analytical continuation.

The same potential gives point convergence. Since $\rev{\Delta_{k+1}}\ge0$,
\begin{equation}\label{eq:lyap-length}
 \ell_{N+j}^2\overset{\eqref{eq:lyap-descent}}{\le}\frac{2\rev{\delta_{N+j}}}L
 \overset{\eqref{eq:lyap-def}}{\le}\frac{2\rev{\Delta_{N+j}}}L
 \overset{\eqref{eq:lyap-rate}}{\le}\frac{2V_N}Lq_{\mathrm{Lyap}}^j.
\end{equation}
The lengths are summable, so the iterates are Cauchy. Their limit belongs to $\X$ by closedness and to $\Xs$ by continuity and \eqref{eq:lyap-rate}. Calling this limit $x_\infty$, the triangle inequality yields
\begin{equation}\label{eq:lyap-tail}
 \|x_{N+n}-x_\infty\|\overset{\eqref{eq:update}}{\le}\sum_{j=n}^{\infty}\ell_{N+j}
 \overset{\eqref{eq:lyap-length}}{\le}
 \frac{\sqrt{2V_N/L}}{1-\sqrt{q_{\mathrm{Lyap}}}}q_{\mathrm{Lyap}}^{n/2}.
\end{equation}
The finite prefix before $N$ does not affect finite length. Since $0<\alpha\le L$, the exponent $-\log q_{\mathrm{Lyap}}$ is of order $\alpha/L$.

\paragraph{Why this proof is shorter, and what is lost.}
The potential combines objective error and squared radius, so there is no need to classify updates or interpolate separate bounds based on their counts. \termrev{The objective decrease bounds the squared step length}, the radius update reduces excess radius, and the invariant ensures that this radius reduction controls a fixed fraction of the potential. The argument still needs the invariant and its initialization proof; it does not shorten those parts of the analysis.

The price is a more conservative geometric estimate: \eqref{eq:lyap-rate} uses the prefactor $V_N=\rev{\Delta_N}+Lt_N^2/2$ and exponent $-\log q_{\mathrm{Lyap}}$, rather than the prefactor and sharper exponent $\lambda$ in \eqref{eq:rlinear-exact}. It also does not recover the separate initialization, contraction, and halving terms in \eqref{eq:complexity}. For example, with $n_\varepsilon=\ceilp{\frac{\log(V_N/\varepsilon)}{-\log q_{\mathrm{Lyap}}}}$, it only gives the combined estimate
\[
 \rev{\Delta_{N+n_\varepsilon}}\auditrev{\overset{\eqref{eq:lyap-rate}}{\le}}
 V_Nq_{\mathrm{Lyap}}^{n_\varepsilon}\le\varepsilon.
\]
This bounds the required additional updates by $n_\varepsilon$. \simplerev{The simpler count above follows from
$1/(-\log q_{\mathrm{Lyap}})\le1/(1-q_{\mathrm{Lyap}})=4(4L+\alpha)/(3\alpha)\le20L/(3\alpha)$, using $\alpha\le L$. Thus $n_\varepsilon\le\widehat n_\varepsilon$.} Thus the Lyapunov argument preserves the rate order and finite-length conclusion, while the counting argument remains preferable for the explicit oracle-call bound and sharper constants.
\endgroup

\subsection{Strong convexity}
{Strong convexity implies quadratic growth, so the preceding guarantees specialize directly to this setting.}

\begin{corollary}[Strongly convex specialization]\label{cor:sc}
\auditrev{Under the smoothness, feasibility, and attainment assumptions of Section~\ref{sec:algorithm},} if $f$ is $\mu$-strongly convex on $\X$ with $0<\mu\le L$, then the minimizer is unique and both \eqref{eq:nonasymptotic} and \eqref{eq:complexity} hold with $\alpha=\mu$, retaining their explicit $t_0$ dependence. For fixed initialization this yields $O\!\left(\frac{L}{\mu}\log\frac{1}{\varepsilon}\right)$ exact-oracle calls, and the iterates converge to the minimizer.
\end{corollary}
\begin{proof}[\normalfont\bfseries Proof of Corollary~\ref{cor:sc}]
A minimizer exists by assumption. It is unique because two distinct minimizers would, by strong convexity, have a midpoint with a strictly smaller objective value. Denote the minimizer by $x_\star$, so $\dist(x,\Xs)=\|x-x_\star\|$.

The constrained first-order condition at $x_\star$ is $\inner{\nabla f(x_\star)}{x-x_\star}\ge0$ for every $x\in\X$. Applying strong convexity at this point gives
\begin{align*}
 f(x)-f_\star
 &\ge\inner{\nabla f(x_\star)}{x-x_\star}+\frac{\mu}{2}\|x-x_\star\|^2\\
 &\ge\frac{\mu}{2}\|x-x_\star\|^2
 =\frac{\mu}{2}\dist(x,\Xs)^2.
\end{align*}
This is quadratic growth with $\alpha_0=\mu$. Since $\mu\le L$, the analysis uses $\alpha=\min\left\{\alpha_0,L\right\}=\mu$. Substituting this value in the two complexity bounds preserves their explicit dependence on $t_0$ and yields the stated asymptotic order. The point-convergence result gives convergence to the unique minimizer.
\end{proof}

\section{Radius variants}\label{app:variants}
\subsection{Reduction to short-step Frank--Wolfe}
\begin{proof}[\normalfont\bfseries Proof of Proposition~\ref{prop:fw}]
For any $x_k,z\in\X$, the diameter bound gives $\|z-x_k\|\le D\le t$. Thus every feasible point lies in the local ball, and
\[
 \X\cap\B(x_k,t)=\X.
\]
The local oracle is therefore exactly the global linear minimization oracle used by \algname{FW}. In particular, the two oracles have the same minimizers and the same gap.

At a nonterminal iteration, $\|z_k-x_k\|>0$. Dividing the accepted length by this positive quantity gives
\[
\begin{relchain}
\frac{\ell_k}{\|z_k-x_k\|} & \overset{\textcolor{black}{\eqref{eq:length}}}{\textcolor{black}{=}} & \frac{1}{\|z_k-x_k\|}
 \min\!\left\{\frac{\inner{\nabla f(x_k)}{x_k-z_k}}{L\|z_k-x_k\|},
                    \|z_k-x_k\|\right\}\\[4pt]
 & = & \min\!\left\{\frac{\inner{\nabla f(x_k)}{x_k-z_k}}{L\|z_k-x_k\|^2},1\right\}.
\end{relchain}
\]
This is the clipped \algname{FW} short-step coefficient. Substituting it into $x_{k+1}=x_k+\frac{\ell_k}{\|z_k-x_k\|}(z_k-x_k)$ gives the \algname{FW} update. Starting from the same point and choosing the same minimizer whenever the oracle is not unique, induction gives identical iterates. The stopping conditions agree because the local and global gaps are equal.
\end{proof}

\subsection{Reduction to projected gradient descent}
\begin{proof}[\normalfont\bfseries Proof of Proposition~\ref{prop:pgd}]
Fix one iteration and omit its index. Let $p=\Proj_\X(x-\gamma\nabla f(x))$ and $t=\|p-x\|$. By definition, $p$ minimizes $\frac{1}{2}\|z-(x-\gamma\nabla f(x))\|^2$ over $z\in\X$. Its constrained first-order condition is
\begin{equation}\label{eq:projcharacter}
 \inner{p-x+\gamma\nabla f(x)}{z-p}\ge0
 \qquad(z\in\X).
\end{equation}
If $t=0$, then $p=x$. Dividing this inequality by $\gamma>0$ gives $\inner{\nabla f(x)}{z-x}\ge0$ for all $z\in\X$, so convexity of $f$ implies that $x$ is optimal.

Suppose $t>0$. The point $p$ belongs to $\X\cap\B(x,t)$. For any other point $z$ in this set, the projection condition yields
\[
\begin{relchain}
\gamma\inner{\nabla f(x)}{z-p} & \overset{\textcolor{black}{\eqref{eq:projcharacter}}}{\textcolor{black}{\ge}} & -\inner{p-x}{z-p}\\[4pt]
 & = & \|p-x\|^2-\inner{p-x}{z-x}\\[4pt]
 & \overset{\textcolor{black}{\eqref{eq:pgdradius}}}{\textcolor{black}{\ge}} & t^2-t\|z-x\|\ge0.
\end{relchain}
\]
The equality uses $z-p=(z-x)-(p-x)$, the next inequality is Cauchy--Schwarz, and the last uses $\|z-x\|\le t$. Dividing by $\gamma>0$ proves that $p$ minimizes the local linear objective.

To prove uniqueness, suppose $z$ is another minimizer. Then the left-hand side of this chain is zero, so all inequalities are equalities. In particular, $\|z-x\|=t$ and $\inner{p-x}{z-x}=t^2$. Consequently,
\[
 \|z-p\|^2=\|z-x\|^2+\|p-x\|^2-2\inner{p-x}{z-x}
 =t^2+t^2-2t^2=0.
\]
Thus $z=p$.

It remains to check that the segment step reaches $p$. Choose $z=x$ in the projection condition and expand it to obtain
\[
 -\|p-x\|^2-\gamma\inner{\nabla f(x)}{p-x}\overset{\textcolor{black}{\eqref{eq:projcharacter}}}{\textcolor{black}{\ge}}0,
 \qquad
 -\inner{\nabla f(x)}{p-x}\overset{\textcolor{black}{\eqref{eq:projcharacter},\eqref{eq:pgdradius}}}{\textcolor{black}{\ge}}\frac{t^2}{\gamma}>0.
\]
The unconstrained minimizing length for the scalar quadratic model is therefore at least
\[
 -\frac{\inner{\nabla f(x)}{p-x}}{Lt}
 \overset{\textcolor{black}{\eqref{eq:projcharacter},\eqref{eq:pgdradius}}}{\textcolor{black}{\ge}}\frac{t}{\gamma L}\ge t,
\]
where the last inequality uses $\gamma\le\frac{1}{L}$. Clipping this length to the segment gives
\[
 \ell\overset{\textcolor{black}{\eqref{eq:length}}}{\textcolor{black}{=}}\min\!\left\{-\frac{\inner{\nabla f(x)}{p-x}}{Lt},t\right\}=t.
\]
Hence the accepted point is $x+\frac{\ell}{t}(p-x)=p$. Repeating this argument at each iteration gives \textcolor{black}{\algname{ProjGD}} with stepsize $\gamma$.
\end{proof}

\subsection{Projected gradient under quadratic growth}\label{app:pgdqg}
\citet{necoara2019} establish linear convergence of \textcolor{black}{\algname{ProjGD}} under quadratic growth. We derive the particular objective contraction used in Table~\ref{tab:comparison} and explain where each assumption enters.

Consider $p=\Proj_\X\left(x-\frac{\nabla f(x)}{L}\right)$. Completing the square in the quadratic upper model gives
\begin{align*}
 &f(x)+\inner{\nabla f(x)}{v-x}+\frac{L}{2}\|v-x\|^2\\
 &\quad=f(x)-\frac{\|\nabla f(x)\|^2}{2L}
       +\frac{L}{2}\left\|v-\left(x-\frac{\nabla f(x)}{L}\right)\right\|^2.
\end{align*}
The first two terms do not depend on $v$, so $p$ minimizes this model over $\X$.

Choose a nearest solution $x_\star$. For $\theta\in[0,1]$, convexity of $\X$ makes $x+\theta(x_\star-x)$ feasible. Smoothness bounds $f(p)$ by the model at $p$, and model optimality bounds this by the model at the comparison point. Thus
\[
\begin{relchain}
f(p) & \rev{\overset{\eqref{eq:smooth}}{\le}} & f(x)+\inner{\nabla f(x)}{p-x}+\frac{L}{2}\|p-x\|^2\\[4pt]
 & \le & f(x)+\theta\inner{\nabla f(x)}{x_\star-x}
             +\frac{L\theta^2}{2}\|x-x_\star\|^2.
\end{relchain}
\]
Convexity of $f$ and quadratic growth give, respectively,
\[
 \inner{\nabla f(x)}{x_\star-x}\le f_\star-f(x),
 \qquad \|x-x_\star\|^2\rev{\overset{\eqref{eq:qg},\eqref{eq:scales}}{\le}}\frac{2(f(x)-f_\star)}{\alpha}.
\]
Substituting these two bounds and collecting terms yields
\[
\begin{relchain}
f(p)-f_\star & \le & (1-\theta)(f(x)-f_\star)
       +\frac{L\theta^2}{2}\frac{2(f(x)-f_\star)}{\alpha}\\[4pt]
 & = & \left(1-\theta+\frac{L\theta^2}{\alpha}\right)(f(x)-f_\star).
\end{relchain}
\]
The choice $\theta=\frac{\alpha}{2L}$ belongs to $[0,1]$ because $0<\alpha\le L$. For this choice,
\[
 1-\theta+\frac{L\theta^2}{\alpha}
 =1-\frac{\alpha}{2L}+\frac{L}{\alpha}\frac{\alpha^2}{4L^2}
 =1-\frac{\alpha}{4L}.
\]
This proves the stated contraction. The comparison point and $\theta$ are used only in the proof. The \textcolor{black}{\algname{ProjGD}} update itself needs $L$, but no growth constant or solution information.

\begin{remark}[Quadratic growth and the squared-gradient-mapping condition]\label{rem:projectedpl}
Define $R_\gamma(x)=\frac{1}{\gamma}(x-\Proj_\X(x-\gamma\nabla f(x)))$. The projected Polyak--\L{}ojasiewicz condition in \citet{locallmo} is $\frac{\|R_\gamma(x)\|^2}{2}\ge\mu_{\mathrm{PL}}(f(x)-f_\star)$ for some $\mu_{\mathrm{PL}}>0$. Quadratic growth does not imply this condition.

To see this, take $\X=[0,1]$, $f(x)=x+\frac{x^2}{2}$, and $L=\gamma=1$. The unique minimizer is zero and $f_\star=0$. Since $x\ge0$,
\[
 f(x)-f_\star=x+\frac{x^2}{2}\ge\frac{x^2}{2}
 =\frac{1}{2}\dist(x,\Xs)^2,
\]
so quadratic growth holds with $\alpha_0=1$. On the other hand, $f'(x)=1+x$, so the gradient step is $x-f'(x)=-1$, whose projection onto $[0,1]$ is zero. Hence $R_1(x)=x$ and, for $x>0$,
\[
 \frac{\|R_1(x)\|^2}{2(f(x)-f_\star)}
 =\frac{x^2}{2x+x^2}=\frac{x}{2+x}\longrightarrow0
 \qquad(x\downarrow0).
\]
The projected Polyak--\L{}ojasiewicz condition would require this ratio to be at least a fixed $\mu_{\mathrm{PL}}>0$, which is impossible near zero. The preceding \textcolor{black}{\algname{ProjGD}} contraction still applies, and in this example the projected step reaches the minimizer immediately.
\end{remark}

\begingroup\color{black}
\subsection{\texorpdfstring{\termrev{Oracle complexity bounds for a common objective accuracy}}{Oracle complexity bounds for a common objective accuracy}}\label{app:tablecurrency}
Table~\ref{tab:comparison} counts calls to each method's stated oracle until the same objective target $\rev{\Delta_k}\le\varepsilon$ is reached. The entries use a common asymptotic convention: problem data and initialization are fixed as $\varepsilon\downarrow0$. Thus finite startup costs and fixed positive scales inside logarithms do not change the displayed order. The oracle types still have different computational costs. We justify the two entries obtained by converting contraction estimates. \rev{The \algname{BroxGD} entry uses the original solution-dependent radius; its inputs footnote records the later addition of backtracking. The displayed total-call bound is not a claim about that later procedure.}

\begin{proof}[\normalfont\bfseries Conversion of the contraction bounds]
Write $\rev{\Delta_k}=f(x_k)-f_\star$. For \algname{ProjGD}, the argument in Appendix~\ref{app:pgdqg} establishes
\begin{equation}\label{eq:table-pgd-step}
 \rev{\Delta_{k+1}}\le\left(1-\frac{\alpha}{4L}\right)\rev{\Delta_k}.
\end{equation}
Consequently,
\begin{equation}\label{eq:table-pgd-error}
 \rev{\Delta_k}\overset{\eqref{eq:table-pgd-step}}{\le}
 \rev{\Delta_0}\left(1-\frac{\alpha}{4L}\right)^k
 \le \rev{\Delta_0}e^{-\alpha k/(4L)}.
\end{equation}
For $0<\varepsilon<\rev{\Delta_0}$, at most $\lceil(4L/\alpha)\log(\rev{\Delta_0}/\varepsilon)\rceil$ updates, and hence projection calls, suffice by \eqref{eq:table-pgd-error}. This gives the displayed $O((L/\alpha)\log(1/\varepsilon))$ bound for fixed initialization.

\rev{For the solution-dependent \algname{BroxGD} row, we use Theorem~3.5 of the original version \citep{locallmo}, not the backtracking theorem of the revised version.} Suppose first that $0<\mu<L$ and set $\auditrev{q_{\mathrm{dist}}}=(L-\mu)/(L+\mu)\in(0,1)$ and $\auditrev{R_0}=\|x_0-x_\star\|$. The cited squared-distance estimate implies
\begin{equation}\label{eq:table-brox-distance}
 \|x_k-x_\star\|\le \auditrev{R_0}\auditrev{q_{\mathrm{dist}}}^k.
\end{equation}
Smoothness at the constrained minimizer, followed by Cauchy--Schwarz, gives
\begin{equation}\label{eq:table-brox-objective}
\begin{relchain}
 \rev{\Delta_k} &\overset{\eqref{eq:smooth}}{\le}&
 \inner{\nabla f(x_\star)}{x_k-x_\star}+\frac L2\|x_k-x_\star\|^2\\[4pt]
 &\le&\|\nabla f(x_\star)\|\,\|x_k-x_\star\|+\frac L2\|x_k-x_\star\|^2\\[4pt]
 &\overset{\eqref{eq:table-brox-distance}}{\le}&
 \|\nabla f(x_\star)\|\auditrev{R_0}\auditrev{q_{\mathrm{dist}}}^k+\frac L2\auditrev{R_0}^2\auditrev{q_{\mathrm{dist}}}^{2k}
 \le C\auditrev{q_{\mathrm{dist}}}^k,
\end{relchain}
\end{equation}
where $C=\|\nabla f(x_\star)\|\auditrev{R_0}+L\auditrev{R_0}^2/2$. If $\auditrev{R_0}=0$, the initial point is already optimal. Otherwise $C>0$. The elementary inequality $-\log(1-u)\ge u$ yields
\[
 -\log\auditrev{q_{\mathrm{dist}}}=-\log\left(1-\frac{2\mu}{L+\mu}\right)
 \ge\frac{2\mu}{L+\mu}\ge\frac{\mu}{L}.
\]
Therefore, for $0<\varepsilon<C$, at most $\lceil(L/\mu)\log(C/\varepsilon)\rceil$ exact local oracle calls suffice by \eqref{eq:table-brox-objective}. If $\mu=L$, the cited distance contraction is zero, so one call suffices from a nonoptimal starting point. This proves the tabulated $O((L/\mu)\log(1/\varepsilon))$ order. The conversion retains the linear term involving $\nabla f(x_\star)$; it does not assume that this gradient vanishes at a constrained optimum.
\end{proof}
\endgroup

\subsection{Convergence with a fixed radius}
\begin{proof}[\normalfont\bfseries Proof of Theorem~\ref{thm:fixed}]
The local oracle and segment step are the same as before. Their existence, feasibility, gap estimates, and objective decrease do not use the adaptive radius update. We first show that the fixed radius eventually satisfies the invariant, then derive the decrease recurrence and the sublinear rate.

While $t<\frac{\alpha}{2L}\rev{s_k}$, the slope estimate from the initialization argument gives
\[
 \frac{\inner{\nabla f(x_k)}{x_k-z_k}}{L\|z_k-x_k\|}
 \overset{\textcolor{black}{\eqref{eq:small-slope}}}{\textcolor{black}{\ge}}\frac{\alpha}{2L}\rev{s_k}
 >t\overset{\textcolor{black}{\eqref{eq:feasible}}}{\textcolor{black}{\ge}}\|z_k-x_k\|.
\]
Thus the step is full. For a boundary step, the same model calculation as in the proof of the invariant reduces $\rev{s_k}$ by at least $\frac{t}{4}$. That calculation uses only the current radius and the accepted point. For an interior step, \termrev{the objective-error estimate for a full interior step} gives
\[
 \frac{\alpha}{2L}\rev{s_{k+1}}
 \overset{\textcolor{black}{\eqref{eq:scales}}}{\textcolor{black}{\le}}\sqrt{\frac{2\rev{\Delta_{k+1}}}{L}}
 \overset{\textcolor{black}{\eqref{eq:reset}}}{\textcolor{black}{\le}}\|z_k-x_k\|<t,
\]
which establishes the invariant at the next point without changing the radius.

If $j$ boundary steps occur before an interior step, summing their decreases gives
\[
 \rev{s_j}
 \overset{\textcolor{black}{\eqref{eq:catchupstep}}}{\textcolor{black}{\le}}\rev{s_0}-\frac{jt}{4}.
\]
The right-hand side is at most $\frac{2Lt}{\alpha}$ when
\[
 j\ge\frac{4}{t}\left(\rev{s_0}-\frac{2Lt}{\alpha}\right).
\]
As in the initialization count, an earlier interior step can only end this phase sooner. Taking the ceiling of the positive part therefore gives at most $C_t$ updates. Once the invariant holds, it persists because $t$ stays fixed and the objective error decreases.

We now derive the recurrence after this phase. First suppose $t\le\rev{s_k}$. Together with the invariant, this is the smaller-radius range in the contraction proof. Its gap and short-step estimates give, for both full and short steps,
\[
 \rev{\Delta_k}-\rev{\Delta_{k+1}}\overset{\textcolor{black}{\eqref{eq:full},\eqref{eq:short},\eqref{eq:slopes},\eqref{eq:small-slope}}}{\textcolor{black}{\ge}}\frac{\alpha}{4L}\rev{\Delta_k}.
\]
In this range $\rev{\Delta_k}\ge\frac{\alpha t^2}{2}$, so
\[
 \frac{\rev{\Delta_k}^2}{2Lt^2}
 \ge\frac{\rev{\Delta_k}}{2Lt^2}\frac{\alpha t^2}{2}
 =\frac{\alpha}{4L}\rev{\Delta_k}.
\]
Thus the linear term is the smaller term in the claimed minimum.

If instead $t>\rev{s_k}$, a nearest solution lies in the local ball. Comparing with it gives a local gap of at least $\rev{\Delta_k}$. A short step therefore satisfies
\[
 \rev{\Delta_k}-\rev{\Delta_{k+1}}
 \overset{\textcolor{black}{\eqref{eq:short}}}{\textcolor{black}{\ge}}\frac{\inner{\nabla f(x_k)}{x_k-z_k}^2}{2L\|z_k-x_k\|^2}
 \overset{\textcolor{black}{\eqref{eq:slopes},\eqref{eq:feasible},\eqref{eq:qg},\eqref{eq:scales}}}{\textcolor{black}{\ge}}\frac{\rev{\Delta_k}^2}{2Lt^2}.
\]
A full step decreases the objective by at least $\frac{1}{2}\rev{\Delta_k}$. Since $\rev{\Delta_k}<\frac{\alpha t^2}{2}\le\frac{Lt^2}{2}$,
\[
 \frac{\rev{\Delta_k}^2}{2Lt^2}
 \overset{\textcolor{black}{\eqref{eq:scales}}}{\textcolor{black}{\le}}\frac{1}{4}\rev{\Delta_k}
 \le\frac{1}{2}\rev{\Delta_k},
\]
so the same quadratic lower bound applies to a full step. Also,
\[
 \frac{\rev{\Delta_k}^2}{2Lt^2}
 \le\frac{\alpha}{4L}\rev{\Delta_k},
\]
so the quadratic term is now the smaller term. Combining the two radius ranges proves the recurrence.

While the error exceeds $\frac{\alpha t^2}{2}$, the recurrence contracts it by the fixed factor $q<1$. It cannot remain above this positive threshold forever, since repeated contraction would bound it by a sequence tending to zero. Hence, unless an optimum is reached earlier, there is a finite index $J$ after initialization with
\[
 0<\rev{\Delta_J}\le\frac{\alpha t^2}{2}.
\]
Monotonicity preserves this upper bound. For all later updates the quadratic term is the smaller one, so
\[
 \rev{\Delta_k}-\rev{\Delta_{k+1}}\overset{\textcolor{black}{\eqref{eq:fixedrec}}}{\textcolor{black}{\ge}}\frac{\rev{\Delta_k}^2}{2Lt^2}.
\]
When both consecutive errors are positive, take their reciprocal difference. Since $\rev{\Delta_{k+1}}\le \rev{\Delta_k}$,
\[
\begin{relchain}
\frac{1}{\rev{\Delta_{k+1}}}-\frac{1}{\rev{\Delta_k}} & = & \frac{\rev{\Delta_k}-\rev{\Delta_{k+1}}}{\rev{\Delta_k}\rev{\Delta_{k+1}}}\\[4pt]
 & \overset{\textcolor{black}{\eqref{eq:descent}}}{\textcolor{black}{\ge}} & \frac{\rev{\Delta_k}-\rev{\Delta_{k+1}}}{\rev{\Delta_k}^2}
 \overset{\textcolor{black}{\eqref{eq:fixedrec}}}{\textcolor{black}{\ge}}\frac{1}{2Lt^2}.
\end{relchain}
\]
Summing from $k=J$ to $J+m-1$ cancels all intermediate reciprocals and gives
\[
 \frac{1}{\rev{\Delta_{J+m}}}
 \overset{\textcolor{black}{\eqref{eq:fixedrec},\eqref{eq:descent}}}{\textcolor{black}{\ge}}\frac{1}{\rev{\Delta_J}}+\frac{m}{2Lt^2}.
\]
Both sides are positive, so taking reciprocals reverses the inequality:
\begin{equation}\label{eq:fixedexplicit}
 \rev{\Delta_{J+m}}\overset{\textcolor{black}{\eqref{eq:fixedrec},\eqref{eq:descent}}}{\textcolor{black}{\le}}
 \left(\frac{1}{\rev{\Delta_J}}+\frac{m}{2Lt^2}\right)^{-1},
 \qquad m\ge0.
\end{equation}
If an optimum is reached after $J$, the bound still holds because the subsequent error is zero. If it is reached before such a $J$ is needed, convergence is already finite.

For $m\ge1$, dropping the positive first term in the denominator gives $\rev{\Delta_{J+m}}\le\frac{2Lt^2}{m}$. With initialization and parameters fixed, $J$ is fixed, and $m=k-J\ge\frac{k}{2}$ for all sufficiently large $k$. Thus $\rev{\Delta_k}=O\!\left(\frac{Lt^2}{k}\right)$ and in particular $f(x_k)\to f_\star$. Finally, quadratic growth gives
\[
 0\le\dist(x_k,\Xs)^2\overset{\textcolor{black}{\eqref{eq:qg},\eqref{eq:scales}}}{\textcolor{black}{\le}}\frac{2\rev{\Delta_k}}{\alpha}\longrightarrow0,
\]
which proves convergence in distance to the solution set.
\end{proof}

\subsection{Transfer of the uniform lower bound}
\begin{proof}[\normalfont\bfseries Proof of Proposition~\ref{prop:lower}]
\citet{halbey2026} provide a fixed strongly convex quadratic $\bar f$ on a Euclidean ball $\bar\X$ and a constant $b_0>0$ with the following property. For every sufficiently large horizon $K$, one can choose a feasible initial point such that short-step \algname{FW} satisfies
\begingroup\makeatletter\def\tagform@#1{\maketag@@@{\textcolor{black}{(#1)}}}\makeatother
\begin{equation}\label{eq:hard-instance}
 \bar f(\bar x_K)-\bar f_\star\ge\frac{b_0}{K^2}.
\end{equation}
\endgroup
The objective and feasible set are fixed, while the initial point may depend on $K$.

Fix $t>0$. Choose $\sigma>0$ such that $\sigma\operatorname{diam}(\bar\X)<t$ and define
\begingroup\makeatletter\def\tagform@#1{\maketag@@@{\textcolor{black}{(#1)}}}\makeatother
\begin{equation}\label{eq:rescaling}
 \X=\sigma\bar\X,\qquad f(x)=\sigma^2\bar f\!\left(\frac{x}{\sigma}\right).
\end{equation}
\endgroup
This scaling is chosen once, independently of $K$, and makes the diameter of $\X$ less than $t$. The chain rule gives
\[
 \nabla f(x)\overset{\textcolor{black}{\eqref{eq:rescaling}}}{\textcolor{black}{=}}\sigma\nabla\bar f\!\left(\frac{x}{\sigma}\right),
 \qquad \nabla^2 f(x)\overset{\textcolor{black}{\eqref{eq:rescaling}}}{\textcolor{black}{=}}\nabla^2\bar f\!\left(\frac{x}{\sigma}\right).
\]
Because $\bar f$ is quadratic, this Hessian is constant. Its eigenvalues, and hence the smoothness and strong-convexity constants, are unchanged. Also, minimizing over the scaled set gives $f_\star=\sigma^2\bar f_\star$.

We next check that scaling preserves the \algname{FW} steps. If $x_k=\sigma\bar x_k$ and $z=\sigma\bar z$, then
\[
 \inner{\nabla f(x_k)}{z}
 \overset{\textcolor{black}{\eqref{eq:rescaling}}}{\textcolor{black}{=}}\sigma^2\inner{\nabla\bar f(\bar x_k)}{\bar z}.
\]
Multiplication by $\sigma^2>0$ preserves minimizers, so corresponding oracle solutions satisfy $z_k=\sigma\bar z_k$. Furthermore,
\[
 \inner{\nabla f(x_k)}{x_k-z_k}
 \overset{\textcolor{black}{\eqref{eq:rescaling}}}{\textcolor{black}{=}}\sigma^2\inner{\nabla\bar f(\bar x_k)}{\bar x_k-\bar z_k},
 \qquad
 \|z_k-x_k\|^2=\sigma^2\|\bar z_k-\bar x_k\|^2.
\]
The factors $\sigma^2$ cancel in the short-step coefficient, while the smoothness constant stays the same. Thus corresponding runs use the same coefficient at each step. From $x_0=\sigma\bar x_0$, induction gives $x_k=\sigma\bar x_k$.

Since the radius $t$ covers $\X$ from every feasible iterate, the fixed-radius method is short-step \algname{FW} by the preceding reduction. For the initialization associated with horizon $K$, it therefore satisfies
\[
 \rev{\Delta_K}
 \overset{\textcolor{black}{\eqref{eq:rescaling}}}{\textcolor{black}{=}}\sigma^2\left(\bar f(\bar x_K)-\bar f_\star\right)
 \overset{\textcolor{black}{\eqref{eq:hard-instance}}}{\textcolor{black}{\ge}}\frac{\sigma^2b_0}{K^2}.
\]
Suppose constants $A<\infty$ and $q_{\mathrm{geom}}\in(0,1)$ bounded the error by $Aq_{\mathrm{geom}}^k$ for every feasible initialization and every $k$. Applying this bound at each of the preceding initializations and multiplying by $K^2$ would give
\[
 0<\sigma^2b_0\le A K^2q_{\mathrm{geom}}^K
\]
for every sufficiently large $K$. The right-hand side tends to zero, because geometric decay dominates the polynomial factor $K^2$. This is a contradiction. The argument rules out a bound uniform over initializations, but does not exclude a geometric bound for an individual initialization with its own prefactor.
\end{proof}

\clearpage
\section{\texorpdfstring{\rev{Frequently used notation}}{Frequently used notation}}\label{app:notation}
\begin{revision}
Table~\ref{tab:notation} collects the notation used most often in the algorithm and its analysis. The objective gap $\Delta_k$ is distinct from the local linear gap $G_k$ and the one-step decrease $\delta_k$. The algorithm does not evaluate $\Delta_k$, $s_k$, or the unknown solution and growth quantities.

\begin{table}[H]
\centering\color{black}
\caption{\rev{Notation at a glance. All norms and projections are Euclidean.}}
\label{tab:notation}
\begingroup
\fontsize{10}{11.2}\selectfont
\setlength{\tabcolsep}{5pt}
\renewcommand{\arraystretch}{1.13}
\newcommand{\notationgroup}[1]{\rowcolor{blue!5}\multicolumn{3}{@{}l@{}}{\strut\bfseries #1}\\\specialrule{0.25pt}{0pt}{3pt}}
\begin{tabular}{@{}>{\raggedright\arraybackslash}p{.19\textwidth}>{\raggedright\arraybackslash}p{.61\textwidth}>{\centering\arraybackslash}p{\dimexpr.20\textwidth-20pt\relax}@{}}
\toprule
\multicolumn{1}{c}{\textbf{Symbol}} & \multicolumn{1}{c}{\textbf{Meaning}} & \textbf{Reference}\\
\midrule
\notationgroup{A. Problem, geometry, and curvature}
$f,\ \X$ & Differentiable convex objective and nonempty closed convex feasible set. & \eqref{eq:oracle}, \eqref{eq:smooth}\\
$f_\star,\ \Xs$ & Minimum objective value and set of minimizers; $\Xs=\argmin_{x\in\X}f(x)$. & Section~\ref{sec:assumptions}\\
$x_\star$ & A minimizer in $\Xs$; uniqueness is not assumed. & Section~\ref{sec:assumptions}\\
\auditrev{$x_{\star,k},\ x_\infty$} & \auditrev{Nearest solution to $x_k$; limit of the adaptive iterates, respectively.} & \auditrev{Lemma~\ref{lem:scaling}, Corollary~\ref{cor:rlinear}}\\
$\B(x,t)$ & Closed Euclidean ball with center $x$ and radius $t$. & \eqref{eq:oracle}\\
$\dist(x,\Xs)$ & Distance from $x$ to the solution set. & Section~\ref{sec:assumptions}\\
$\Proj_C(v)$ & Euclidean projection of $v$ onto a nonempty closed convex set $C$. & \eqref{eq:segment}\\
$L$ & Known positive constant in the quadratic upper model for $f$. & \eqref{eq:smooth}\\
$\alpha_0,\ \alpha$ & \termrev{Quadratic-growth constant and its truncation at $L$: $\alpha=\min\left\{\alpha_0,L\right\}$.} & \eqref{eq:qg}, \eqref{eq:scales}\\
$\mu$ & Strong-convexity constant, when strong convexity is assumed. & Table~\ref{tab:comparison}\\
\midrule
\notationgroup{B. Quantities in one algorithm update}
$x_k,\ t_k$ & Current feasible iterate and positive local-oracle radius. & Algorithm~\ref{alg:main}\\
$g_k,\ z_k$ & Gradient $g_k=\nabla f(x_k)$ and an exact linear minimizer over $\X\cap\B(x_k,t_k)$. & \eqref{eq:oracle}, \eqref{eq:stepnotation}\\
$d_k,\ D_k$ & Oracle displacement $d_k=z_k-x_k$ and its length $D_k=\|d_k\|$. & \eqref{eq:stepnotation}\\
$G_k$ & Local linear gap $-\inner{g_k}{d_k}$; zero certifies optimality. & \eqref{eq:stepnotation}\\
$\eta_k$ & Clipped step coefficient $\min\left\{1,\frac{G_k}{LD_k^2}\right\}$ at a nonterminal update. & \eqref{eq:coefficient}\\
$\ell_k$ & Accepted displacement length $\|x_{k+1}-x_k\|=\eta_kD_k$. & \eqref{eq:length}\\
$S_k$ & Segment $[x_k,z_k]$ used in the projection interpretation. & \eqref{eq:segment}\\
\midrule
\notationgroup{C. Objective accuracy and convergence analysis}
$\Delta_k$ & Objective gap $f(x_k)-f_\star$. & \eqref{eq:analysisnotation}\\
$s_k$ & \textcolor{black}{Upper bound on solution distance} $\sqrt{2\Delta_k/\alpha}$; $\dist(x_k,\Xs)\le s_k$. & \eqref{eq:analysisnotation}\\
$\delta_k$ & One-step objective decrease $\Delta_k-\Delta_{k+1}$. & Appendix~\ref{app:invariant}\\
$\varepsilon,\ K(\varepsilon)$ & Positive objective-error target and first index with $\Delta_{K(\varepsilon)}\le\varepsilon$. & Theorem~\ref{thm:main}\\
$N,\ C_0$ & An index establishing $t_N\ge\alpha s_N/(2L)$ and the initialization bound $N\le C_0$ for a suitable such index. & \eqref{eq:invariant}, \eqref{eq:catchup}\\
$q$ & Objective contraction factor $1-\alpha/(4L)$. & \eqref{eq:q}\\
\auditrev{$R_0,\ q_{\mathrm{dist}}$} & \auditrev{Original \algnametiny{BroxGD} bound: initial solution distance and distance-contraction factor, respectively.} & \auditrev{Appendix~\ref{app:tablecurrency}}\\
$P_\varepsilon,\ H_\varepsilon$ & \termrev{Bounds on updates satisfying \eqref{eq:contract} and on remaining updates (radius halvings) before reaching the target.} & \eqref{eq:productivecount}, \eqref{eq:halvingcount}\\
$r_k,\ V_k$ & \termrev{Scaled squared radius $r_k=Lt_k^2/2$; Lyapunov function $V_k=\Delta_k+r_k$.} & \eqref{eq:lyap-def}\\
$\log,\ \log_2$ & Natural logarithm ($\log=\ln$) and base-two logarithm, respectively. & \eqref{eq:initialcomplexity}\\
$[v]_+$ & Positive part $\max\left\{v,0\right\}$. & Section~\ref{sec:convergence}\\
\bottomrule
\end{tabular}
\endgroup
\end{table}
\end{revision}

\end{document}